\documentclass{article}                    
\usepackage{graphicx}
\usepackage{latexsym}
\usepackage{amssymb}
\usepackage{theorem}
\usepackage{graphicx}
\newcommand{\R}{\hbox{I}\!\hbox{R}}
\newtheorem{teo}{Theorem}[section]
\newtheorem{cor}{Corollary}[section]
\newtheorem{lem}{Lemma}[section]
\newtheorem{prop}{Proposition}[section]
\newtheorem{pro}{Proof}
\newtheorem{conj}{Conjecture}

\title{On uniqueness of the $q$-state Potts model on a self-dual family of graphs}

\author{J.-M. Billiot \thanks{LJK, UMR 5224, BSHM, Universit\'e Pierre Mend\`es France, 1251 Avenue Centrale, BP 47, 38040 Grenoble Cedex 9, France}
 \and F. Corset \and E. Fontenas
}

\begin{document}
\maketitle

\begin{abstract}
This paper deals with the location of the complex zeros of the Tutte polynomial for a class of self-dual graphs. 
For this class of graphs, as the form of the eigenvalues is known, 
the regions of the complex plane can be focused on the sets where there is only one dominant eigenvalue in particular 
containing the positive half plane. Thus, in these regions, the analyticity of the pressure can be derived easily. 
Next, some examples of graphs with their Tutte polynomial having a few number of eigenvalues are given. 
The cases of the strip of triangles with a double edge, 
the wheel and the cycle with an edge having a high order of multiplicity are presented. 
In particular, for this last example, we remark that the well known conjecture of Chen et al. \cite{Chen96} is false in the finite case.
\end{abstract}

\section{Introduction}
The program of studying complex zeros of the partition function was pioneered by Yang and Lee \cite{Lee52}. 
The location of the complex zeros of the Tutte polynomial is important from a statistical mechanic point of view because it is related to 
a possible phase transition for q-Potts model on this family of graphs using the Fortuin Kasteleyn representation 
(for some correlation duality relations for the planar Potts model, see \cite{King02a}). 
There exists an impressive literature concerning the location of zeros of the Tutte and chromatic polynomials for a large class of graphs and, 
in this context, the Beraha numbers play an important role (see for example \cite{Salas01} and references therein). 
Sokal \cite{Sokal04} proves for general graphs that these zeros are dense in the complex plane. 
A large number of conjectures for a wide class of graphs such as triangulation, 
planar, cubic  graphs are given in \cite{Jackson03}. For self dual graphs, an intriguing conjecture proposed by Chen and al. 
\cite{Chen96} asserts that, in a half plane, the complex zeros are located on a circle. 

First, our goal was to prove this conjecture for the large class of self dual graphs as possible. In fact, in the finite case, we believed that the self duality, implying the symmetry of the Tutte polynomial, was the adequate property in order to obtain such a result. More precisely, we thought that it forced in the positive half plane the complex zeros of this polynomial to be located on the unit circle in the same spirit as the Lee and Yang theorem. We realized that it is not true. We give an example and surprisingly we find a sequence of self dual graphs for which this unit circle does not belong to the accumulation set of the zeros. But the question raised by this conjecture remains at the thermodynamic limit. Unfortunately in the general case, for example for self dual strips of the square lattice, the eigenvalues coming from the powerful transfer matrix method are roots of polynomials with high degree: then, it is difficult to study the location of curves of degeneration of the dominant eigenvalue even if we focus on the positive half plane.

In this work, we choose to study a family of self dual graphs. For this class of graphs, as the form of the eigenvalues is known, we are able to present regions having only one dominant eigenvalue in particular containing the positive half plane.
Then we derive easily in these regions of the complex plane the analyticity of the pressure. Next, we study some examples containing a few number of eigenvalues. We use deletion and contraction rules of the Tutte polynomial to obtain recursive formula rather than the transfer-matrix method. One originality of this approach is to unify in the same framework different kind of self dual corrections. Besides, we use a suitable variable allowing us to identify the set of accumulation of zeros more easily.
 
The outline of this paper is the following. In the first section, the conjecture of \cite{Chen96} and the well-known results such as Beraha Kahane Weiss theorem and Vitali's convergence theorem are recalled. Next, general results, giving regions where we have only one dominant eigenvalue, are established and the analyticity of the pressure is deduced. Then for the strip of triangles with a double edge, the wheel, the cycle with an edge having a high order of multiplicity, the set of degeneration of the dominant eigenvalue is obtained. To conclude, some perspectives and conjectures raised by ours results are proposed. The technical aspects of the proofs are given in the last section.

\section {Preliminaries and framework}

In this paper, we deal with self dual graphs ($G=G^{\star}$) where $G^{\star}$ denotes the dual graph of $G$. One useful tool in graph theory is the Tutte or dichromatic polynomial. It is a function of two variables $x$ and $y$; we denote $T(G,x,y)$ the Tutte polynomial of the graph $G$ calculated at the point $(x,y)$.
A well known fact  is that for self dual graphs, this function is symmetric $$T(G,x,y)=T(G^{\star},x,y)=T(G,y,x).$$ On a statistical mechanic point of view, the link with Ising and $q$-Potts models using the Fortuin Kasteleyn representation leads to look in particular on the hyperbola $(x-1)(y-1)=q$.  In all the following, we assume that $q\geq 1$.  Usually, we consider the variable $v=re^{i\theta}$ by taking $$x=1 + \frac{\sqrt{q}}{v}\quad\mbox{and} \quad y=1 + \sqrt{q}{v}.$$
In our case, we choose working with the variable $$z=x+y-2=\sqrt{q}({v}+1/v)\quad \mbox{then}\quad xy=z+q+1.$$ It is easy to see that, on this curve and from the symmetry, the Tutte polynomial is also a polynomial of the variable $z$. Chen and al. \cite{Chen96} propose the following intriguing conjecture:
\begin{conj}
For finite planar self-dual lattices and for square lattice with free or periodic boundary conditions in the thermodynamic limit, the Potts partition zeros in the $Re(v)>0$ half plane are located on the unit circle $|v|=1$.
\end{conj}

\begin{figure}[ht]
\begin{center}
\includegraphics[width=4cm]{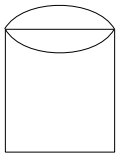}
\caption{Graph $C$ : a cycle with an edge of multiplicity three}
\label{grapheC}
\end{center} 
\end{figure}

Our first idea was to prove this conjecture on the finite case on the large class of self dual graphs as possible. But later we realize that this conjecture is not true in the finite case. Let us take the following simple example: choose the self dual graph $C$ with four vertices and six edges defined as a cycle of length  four and an edge of multiplicity three, see Figure \ref{grapheC}.
The Tutte polynomial of this graph on the previous hyperbola with $q=16$ can be written as a function of $z$ as
$$T(C,z)=z^3+6z^{2}+12z+218.$$

There are one not positive real root ($-2-(210)^{1/3}$) and two conjugated roots with not negative real part ($(210)^{1/3}/2-2\pm i(210)^{1/3} \sqrt{3}/2$); then, from the relation between $z$ and $v$, the complex roots in $z$ with $Re(z)>0$ correspond with the roots in $v$ in the region $Re(v)>0$ not on the unit circle. At least in the finite case, this conjecture is false. \\

In this work, at the thermodynamic limit, the analyticity of the pressure can be obtain by studying and avoiding the location of the accumulation sets of zeros of the partition function (or Tutte polynomial) for a family of self-dual graphs.\\

A central role in our work is played by a theorem on analytic functions due to
Beraha, Kahane and Weiss. Let $D$ be a domain (connected open set) in the complex plane,
and let $\alpha_{1}, \ldots, \alpha_{M}, \lambda_{1},\ldots,\lambda_{M}$ be analytic functions on D, none of which is
identically zero. For each integer $n\geq0$, let define
$$f_{n}(z) =
\sum _{k=1}^{M}
\alpha_{k}(z)[\lambda_{k}(z)]^n.$$
We are interested in the zero sets
$ Z(f_{n}) = \{z\in D: f_{n}(z) = 0\}$
and in particular in their limit sets as $n \rightarrow\infty$:\\

\noindent - liminf $Z(f_{n})=\{z\in D:$ every neighborhood $U\ni z$ has a nonempty intersection
with all but finitely many of the sets $Z(f_{n})\}$.\\

\noindent - limsup $Z(f_{n})=\{z \in D:$ every neighborhood $U\ni z$ has a nonempty intersection
with infinitely many of the sets $Z(f_{n})\}$.\\

Let $k$ be a dominant subscript z if $|\lambda_{k}(z)| \geq |\lambda_{l}(z)|$ for all $l\in\{1 \ldots M\}$.
Then the limiting zero sets can be completely characterized as follows:
\begin{teo}
Let $D$ be a domain in $\mathbb C$, and let $\lambda_{1},\ldots \lambda_{M}, \alpha_{1},\ldots, \alpha_{M}$  be analytic functions on $D$, none of which is identically zero. Let us further
assume a no-degenerate-dominance condition: there do not exist subscripts $k \neq k^{'}$
such that $\lambda_{k}=\omega \lambda_{k^{'}}$ for some constant $\omega$ with $|\omega| = 1$ and such that $\{z \in D: k\, \mbox{is dominant at}\, z\}$ $( = \{z \in D: k'\, \mbox{is dominant at}\, z\})$ has
nonempty interior. For each integer $n \geq0$, define $f_n$ by
$$f_{n}(z) =
\sum _{k=1}^{M}
\alpha_{k}(z)\lambda_{k}(z)^{n} $$

Then $\mbox{\it{liminf}}\, Z(f_{n}) = \mbox{\it{limsup}}\, Z(f_{n})$, and a point $z$ lies in this set if and only if either\newline
(a) There is a unique dominant subscript $k$ at $z$, and $\alpha_{k}(z) = 0$; or\newline
(b) There are two or more dominant subscripts at $z$.
\end{teo}

Note that case (a) consists of isolated points in $D$, while case (b) consists of curves
(plus possibly isolated points where all the $\lambda_k$ vanish simultaneously). 

We also recall the classical famous Vitali's convergence theorem useful to obtain analyticity properties:
\begin{teo}
Let $p_{n}(z)$ be a sequence of functions, each regular in a region $D$. Assume that it exists a constant $B$ as $|p_{n}(z)|\leq B$ for every $n$ and  for all $z\in D$. If $p_{n}(z)$ tends to a limit as $n\rightarrow\infty$ at a set of points having a limit point inside $D$, then $p_{n}(z)$ tends uniformly to a limit in any region bounded by a contour interior to $D$: the limit therefore being an analytic function~of~$z$. 
\end{teo}
 
\section{General results}
Here, we are interested in the following family:
$$f_{n}(z) =
\sum _{k=1}^{M}\Bigl[
\alpha_{a_k}(z)[\lambda^{+}_{a_k}(z)]^{n}+ \beta_{a_k}(z)[\lambda^{-}_{a_k}(z)]^{n}\Bigr]$$
where $\lambda^{+}_{a_k}(z)$ and $\lambda^{-}_{a_k}(z)$ are the solutions of the following equation 
$$X^{2}-(z+2+a_{k})X+ z+q+1=0$$ with $a_{k}\in[0,q]$. They can be expressed as follows
$$\lambda_{a_k}^{\pm}(z)=\frac12\left(z+2+a_{k}\pm\displaystyle\sqrt{(z+a_{k})^2-4(q-a_{k})}\right).$$
By convention, we denote $\lambda_{a_k}(z)$, the solution between $\lambda^{+}_{a_k}(z)$ and $\lambda^{-}_{a_k}(z)$ with the greatest magnitude. The functions $\{\alpha_{a_k},\beta_{a_k},\,k=1..M\}$ are chosen such that the function $f_n$ stays a polynomial function in the variable $z$.
We denote by
$a_{u}=\sup_{k=1\ldots M}a_{k}$ and $a_{l}=\inf_{k=1\ldots M}a_{k}$. We also use some subsets of the complex plane: $$D_{]-\infty,-a_{u}-2[}=\{z=c+id; (c,d)\in\R^2,\, c<-a_{u}-2\},$$
$$ D_{]-a_l,+\infty[}=\{z=c+id; (c,d)\in\R^2,\, c>-a_{l}\}.$$

\noindent Then, the analyticity of the pressure in some regions can be derived easily:

\begin{teo}\label{theo1}
There exists only one dominant eigenvalue at $z$:\\
- $\forall z\in D_{]-a_l,+\infty[} \setminus\{c\in[-a_l, \sup(-a_l, -a_{u}+2\sqrt{q-a_{u}})], \,d=0\}$, $\lambda_{a_u}(z)$ is the dominant eigenvalue.\\
- $\forall z\in D_{]-\infty,-a_{u}-2[}\setminus\{c\in[\inf(-a_u-2, -a_{l}-2\sqrt{q-a_{l}}),-a_u-2], \,d=0\}$, $\lambda_{a_l}(z)$ is the dominant eigenvalue. 
\end{teo}
\noindent{\bf Proof:}
\rm It is a direct consequence of lemma \ref{lem1} and \ref{lem2}.\\

\noindent Now taking $p_{n}(z)=\displaystyle\frac{\ln(f_{n}(z))}{n}$, we can deduce the following corollary.

\begin{cor} We have
\begin{itemize}
\label{coro1}
\item If $\forall z\in D_{]-a_l,+\infty[} \setminus\{c\in[-a_l, \sup(-a_l, -a_{u}+2\sqrt{q-a_{u}})], \,d=0\}$, $\alpha_{a_u}(z)\neq 0$ then $p_{n}(z)\longrightarrow \ln[\lambda_{a_u}(z)]$ as $n\rightarrow\infty$.
\item If $\forall z\in D_{]-\infty,-a_{u}-2[}\setminus\{c\in[\inf(-a_u-2, -a_{l}-2\sqrt{q-a_{l}}),-a_u-2], \,d=0\}$, $\beta_{a_l}(z)\neq 0$ then $p_{n}(z)\longrightarrow \ln[\lambda_{a_l}(z))]$ as $n\rightarrow\infty$.
\end{itemize}

Let $K$ be a non negative real and $B((0,0),K)$ the disk center at $(0,0)$ with radius $K$. Both limits are analytic functions of $z$ respectively on subsets $D_{]-a_l,+\infty[} \setminus\{c\in[-a_l, \sup(-a_l, -a_{u}+2\sqrt{q-a_{u}})], \,d=0\}\bigcap B((0,0),K)$ and $D_{]-\infty,-a_{u}-2[}\setminus\{c\in[\inf(-a_u-2, -a_{l}-2\sqrt{q-a_{l}}),-a_u-2], \,d=0\}\bigcap B((0,0),K)$.
\end{cor}
\noindent{\bf Proof:}
From the previous theorem, $\forall z\in D_{]-a_l,+\infty[} \setminus\{c\in[-a_l, \sup(-a_l, -a_{u}+2\sqrt{q-a_{u}})], \,d=0\}$, the dominant eigenvalue is $\lambda_{a_u}(z)$. Because of the given form of $\lambda_{a_u}(z)$ we obtain that$$1\leq|\lambda_{a_u}(z)|\leq C_{1}(K)$$ with $C_{1}(K)$ is a constant depending only of $K, a_{u}$ and $q$. We have also that $\forall z\in D_{]-\infty,-a_{u}-2[}\setminus\{c\in[\inf(-a_u-2, -a_{l}-2\sqrt{q-a_{l}}),-a_u-2], \,d=0\}$, the dominant eigenvalue is $\lambda_{a_l}$. We find in this region that $$\frac{a_{u}-a_{l}}{2}\leq|\lambda_{a_l}(z)|\leq C_{2}(K)$$ with $C_{2}(K)$ is a constant depending only of $K, a_{l}$ and $q$. Now the result is a direct consequence of Vitali's convergence theorem.
\section{Examples}
After working in two specified regions of the complex $z$ plane, it may be be interesting to know what happens elsewhere. It is the goal of this section. Several simple examples of self-dual graphs are given: the classical one is the wheel. The other graphs - the strip of triangles with a double edge, the cycle with one multiple edge - are built from the wheel by moving some of its edges. By making this transformation, we only have to keep the symmetry property of the Tutte polynomial. We provide their Tutte polynomial using the contraction and deletion rules useful to obtain recurrence formula. We need these recurrence formula to identify without ambiguity the functions $\lambda_{a_k}^{\pm}$ introduced before. These graphs belong to the framework we discussed before for a particular choice of the parameters $a_k$ like 0, 1 or $q$. We present the location of accumulation sets of zeros and more precisely the curves describing the degeneration of the dominant eigenvalue using respectively the variable $z$ and $v$ in the whole complex plane. For all graphs $G$, the notation $T(G)$ for the Tutte polynomial $T(G,x,y)$ evaluated at the point $(x,y)$ shall be used.  
\subsection{Triangles with a double edge}
\begin{figure}[ht]
\begin{center}
\includegraphics[width=6cm]{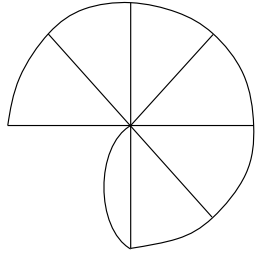}
\caption{$G_n$ a strip of triangles with a double edge}\label{triangle}
\end{center} 
\end{figure}

These graphs are called $G_{n}$ (see Figure \ref{triangle}) and graphs having $n$ triangles $Tr_{n}$. 
We have these relations:
$$\left\{\begin{array}{l}
T(G_n)=T(Tr_n)+yT(G_{n-1})\\
T(Tr_n)=xT(Tr_{n-1})+T(G_{n-1})
\end{array}\right.\Longleftrightarrow
\left\{\begin{array}{l}
T(G_n)=T(Tr_n)+yT(G_{n-1})\\
yT(Tr_n)=xyT(Tr_{n-1})+yT(G_{n-1})
\end{array}\right.$$
This implies that
\begin{equation}\label{eqgn}
T(G_n)=(1+y)T(Tr_n)-xyT(Tr_{n-1}).
\end{equation}
We prove easily the following recurrence formula for Tutte polynomial of graphs $Tr_n$
\begin{equation}\label{recTr}
T(Tr_n)=(1+x+y)T(Tr_{n-1})-xyT(Tr_{n-2}).
\end{equation}
Using (\ref{eqgn}) and (\ref{recTr}), it comes
\begin{equation}\label{recGn}
T(G_n)=(1+x+y)T(G_{n-1})-xyT(G_{n-2}).
\end{equation}
We deduce that:
$$U_n=\left(\begin{array}{c}T(G_n)\\T(G_ {n-1})\end{array}\right)=\left(\begin{array}{cc}1+x+y&-xy\\1&0\\\end{array}\right)\left(\begin{array}{c}T(G_{n-1})\\T(G_ {n-2})\end{array}\right).$$
The eigenvalues of this matrix are $$\mu_1=\frac{x+y+1+\sqrt{(1+x+y)^2-4xy}}2\quad \mu_2=\frac{x+y+1-\sqrt{(1+x+y)^2-4xy}}2.$$
The matrix defining eigenvectors is denoted by $P$ and its inverse $P^{-1}$. They are given by :
$$P=\left(\begin{array}{cc}\mu_1&\mu_2\\1&1\end{array}\right),\quad P^{-1}=\frac{1}{\mu_1-\mu_2}
\left(\begin{array}{cc}1&-\mu_2\\-1&\mu_1\end{array}\right).$$
Now, it comes
$$U_n=\left(\begin{array}{cc}\mu_1&\mu_2\\1&1\end{array}\right)\left(\begin{array}{cc}\mu_1^{n-1}&0\\0&\mu_2^{n-1}\\\end{array}\right)
\frac{1}{\mu_1-\mu_2}\left(\begin{array}{cc}1&-\mu_2\\-1&\mu_1\end{array}\right)U_1$$ where
$$U_1=\left(\begin{array}{c}\mu_1^2+\mu_2^2+\mu_1\mu_2-\mu_1-\mu_2\\\mu_1+\mu_2-1\end{array}\right).$$
Then, the Tutte polynomial associated with $G_n$ is
 $$T(G_n)=\frac{\mu_1-1}{\mu_1-\mu_2}\mu_1^{n+1}+\frac{\mu_2-1}
{\mu_2-\mu_1}\mu_2^{n+1}.$$
\noindent On the hyperbola $(x-1)(y-1)=q$, these eigenvalues are of the form introduced before with $a=1$ and the Tutte polynomial can be written as $$T(G_n)=\alpha_{1}(z)[\lambda^{+}_{1}(z)]^{n+1}+\beta_{1}(z)[\lambda^{-}_{1}(z)]^{n+1}$$
with $\alpha_{1}(z)=\frac{\lambda^{+}_{1}(z)-1}{\lambda^{+}_{1}(z)-\lambda^{-}_{1}(z)}$ and $\beta_{1}(z)=\frac{\lambda^{-}_{1}(z)-1}{\lambda^{-}_{1}(z)-\lambda^{+}_{1}(z)}$.\newline
 It is outstanding that such a simple form of Tutte polynomial leads to a large variety of cases according to the values of $q$ when studying the degeneration of the dominant eigenvalue. More precisely, we have both following propositions respectively using complex variable $z$ and $v$.
\begin{prop}
For the family of graphs $(G_{n})_{n\geq 0}$, the location of the degeneration of the dominant eigenvalue is described in the complex plane using variable $z=c+i\,d$ as follows:\\
- $\forall\, q\geq 2$
$$|\lambda^{+}_{1}(z)|=|\lambda^{-}_{1}(z)| \Longleftrightarrow d=0,\, c\in[-1-2\sqrt{q-1},-1+2\sqrt{q-1}]$$

\noindent - $\forall\, q\in[1,2]$
$$|\lambda^{+}_{1}(z)|=|\lambda^{-}_{1}(z)| \Longleftrightarrow \left\{\begin{array}{l}d=0,\, c\in[-1-2\sqrt{q-1},-1+2\sqrt{q-1}]\\
\mbox{or}\,\, z\in C((-q-1,0),2-q)\end{array}\right.$$
where $C((-q-1,0),2-q)$ denoted the circle of center $(-q-1,0)$ and of radius $2-q$.
\end{prop}

\noindent{\bf Proof: }we just have to apply lemma \ref{lem2} in the particular case $a=1$.\\

\noindent By using lemma \ref{lem4} and lemma \ref{lem5}, it comes

\begin{figure}[ht!]
\begin{center}
\begin{tabular}{cc}
\includegraphics[width=4cm,angle=-90]{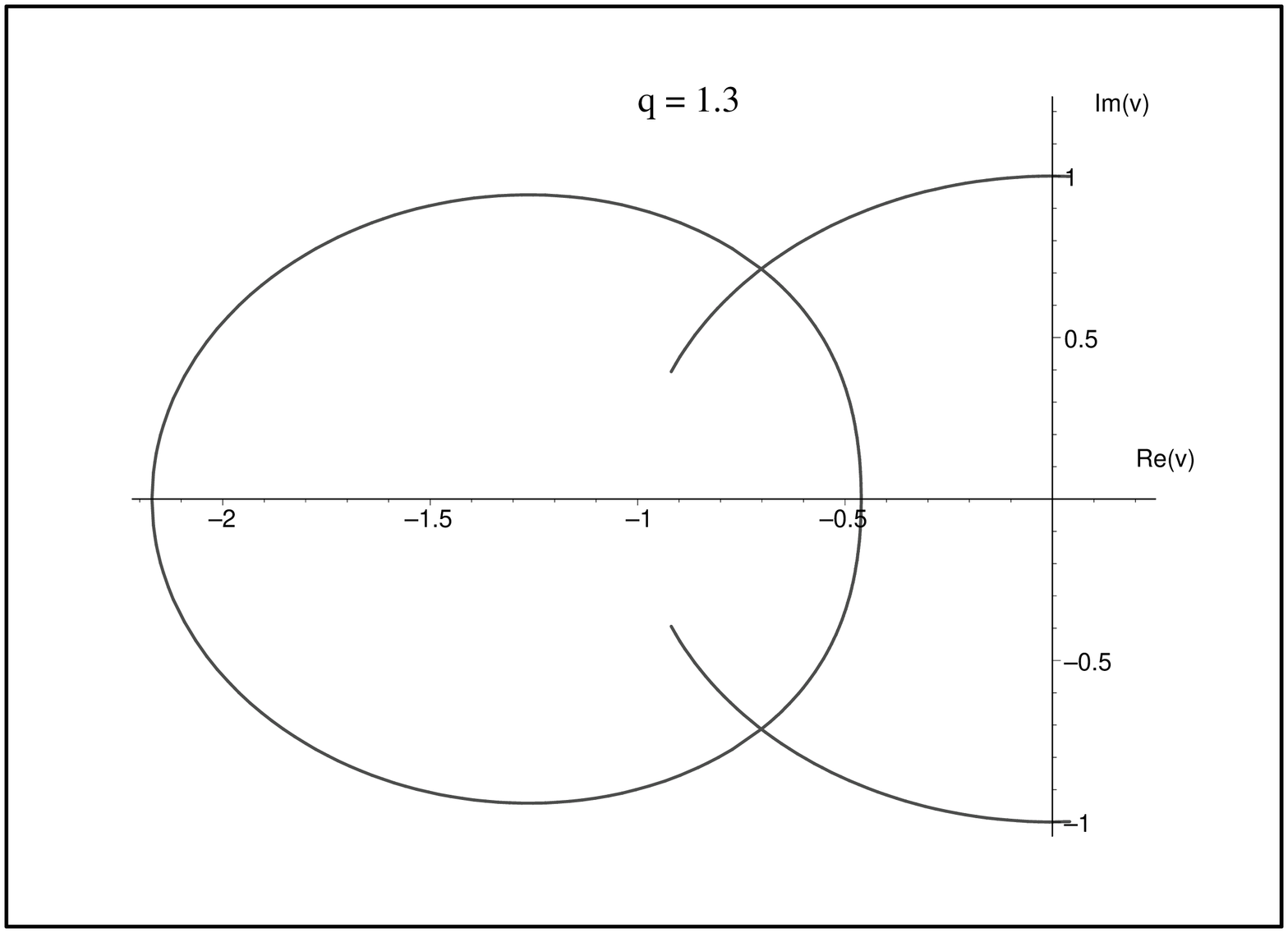}&\includegraphics[width=4cm,angle=-90]{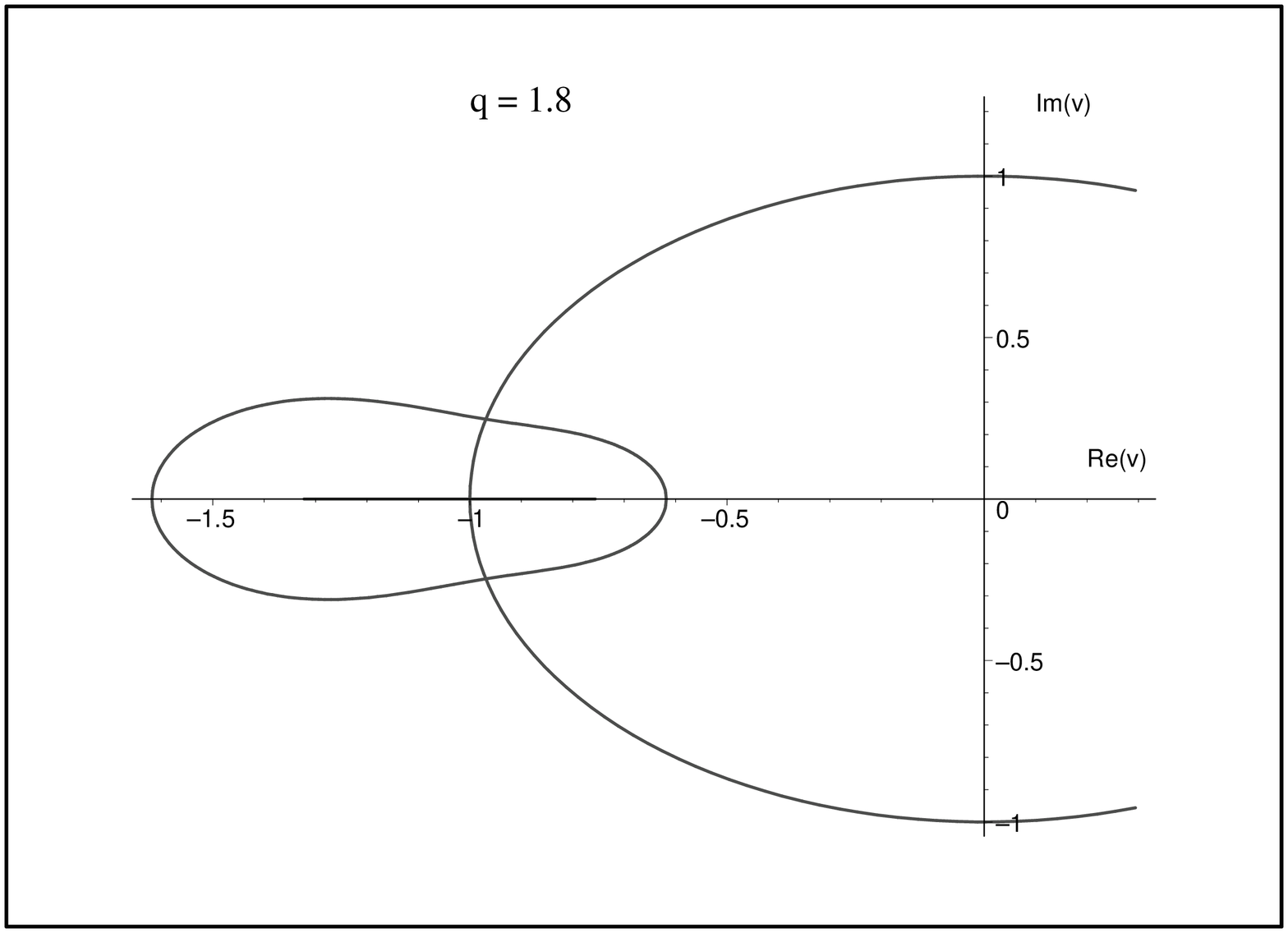}\\
\includegraphics[width=4cm,angle=-90]{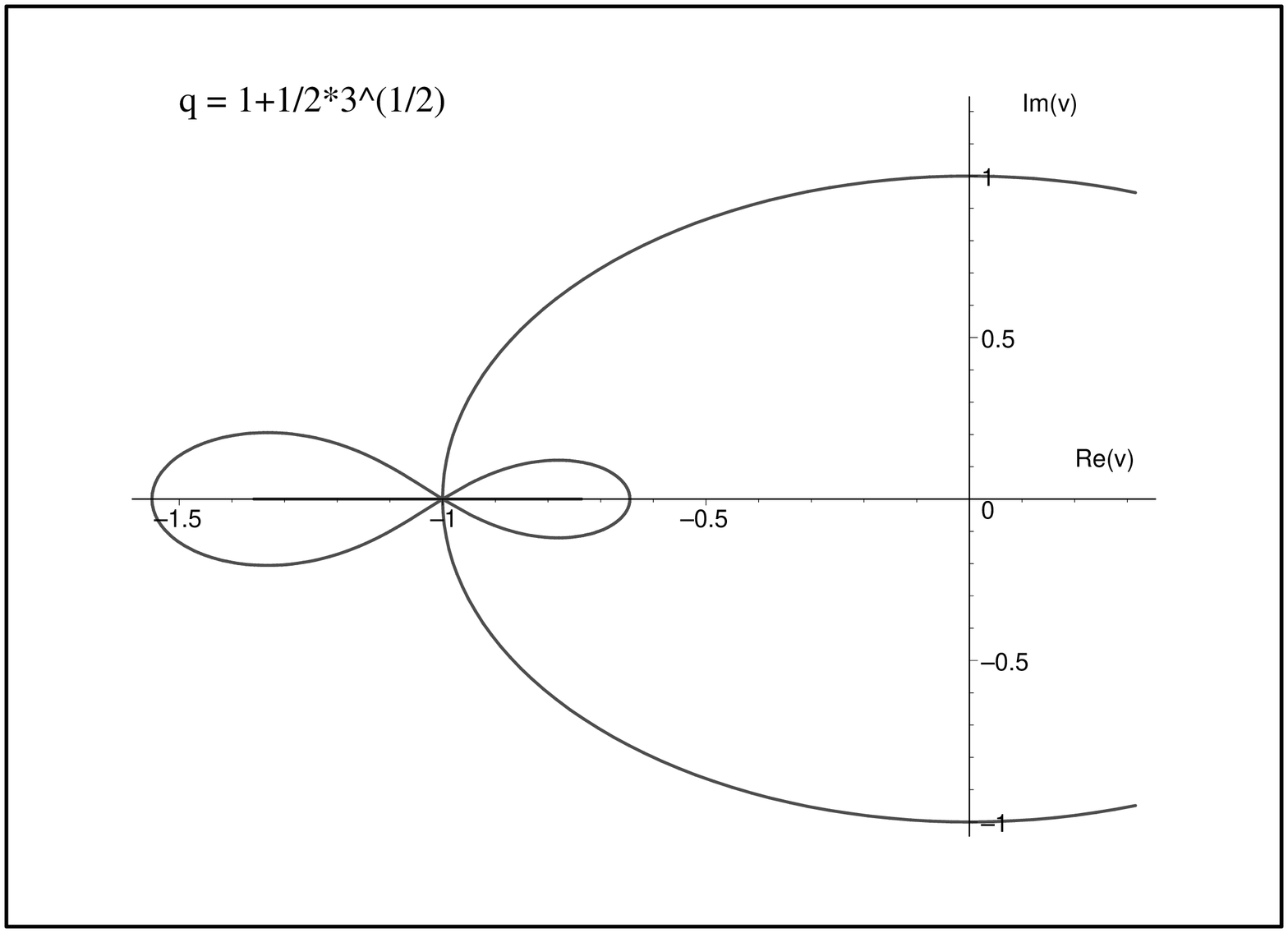}&\includegraphics[width=4cm,angle=-90]{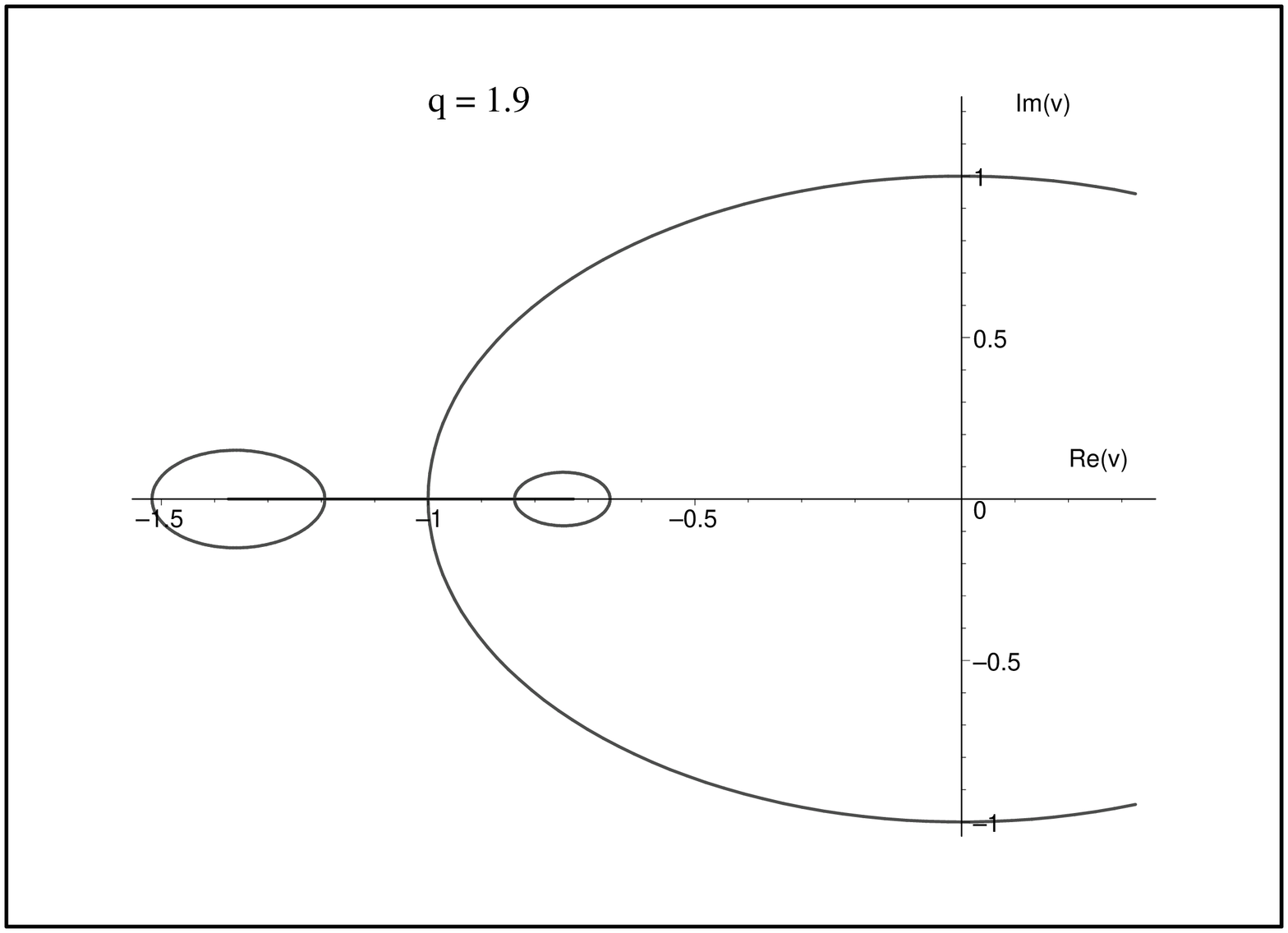}
\end{tabular}
\end{center}
\caption{The limiting zero sets for distinct values of $q$ in the $v$ complex plane.} 
\label{qGn}
\end{figure}

\begin{prop}
For the family of graphs $(G_{n})_{n\geq 0}$, the location of the degeneration of the dominant eigenvalue is described in the complex plane using the variable $v=re^{i\theta}$ as follows:\\
We have $|\lambda^{+}_{1}(z)|=|\lambda^{-}_{1}(z)|$ when:\\

\noindent- if $q\in[1,25/16]$ $\left\{\begin{array}{l}\theta\in[ \arccos(\frac{-1+2\sqrt{q-1}}{2\sqrt{q}}),\arccos(\frac{-1-2\sqrt{q-1}}{2\sqrt{q}})],\,r=1\\
\mbox{or}\,\, v\in F(C((-q-1,0),2-q)).\end{array}\right.$\\

\noindent- if $q\in[25/16,2]$  $\left\{\begin{array}{l}\theta\in[ \arccos(\frac{-1+2\sqrt{q-1}}{2\sqrt{q}}), \pi],\, r=1\\
\mbox{or}\,\,\theta=\pi,\,r\in[1/r_{1},r_{1}]\\
\mbox{or}\,\, v\in F(C((-q-1,0),2-q)).\end{array}\right.$\\

\noindent- if $q\geq 2$ $\left\{\begin{array}{l}\theta\in[\arccos(\frac{-1+2\sqrt{q-1}}{2\sqrt{q}}), \pi],\, r=1\\
\mbox{or}\,\,\theta=\pi,\,r\in[1/r_1;r_1]\end{array}\right.$\\
with $r_{1}$ and $1/r_1$ the roots of the polynomial $\sqrt{q}r^{2}-(1+2\sqrt{q-1})r+\sqrt{q}=0$.
\end{prop}
The Figure \ref{qGn} shows  for this family of graphs, the location of the degeneration of the dominant eigenvalue in the $v$ complex plane for different values of $q$. We can notice the particular value $q_0=1+\sqrt{3}/2$.

\subsection{Wheel graphs}

\begin{figure}[ht]
\begin{center}
\includegraphics[width=6cm]{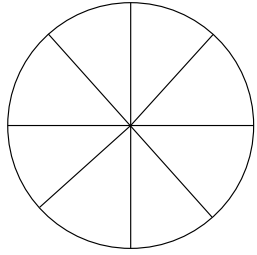}
\caption{The wheel}
\label{wheel}
\end{center} 
\end{figure}
We consider the wheel, a planar graph with $n+1$ vertices with a cycle of length $n$ and one internal vertex linked to all the vertices of the cycle (see Figure \ref{wheel}). Then, we obtain a wheel graph  having $n$ triangles, denoted $B_n$. The Tutte polynomial for these graphs can be deduced from the relations between the Tutte polynomials of $G_n$ and $Tr_n$. We have
$$T(B_n)=T(Tr_{n-1})+T(B_{n-1})+y(T(G_{n-2})-xT(G_{n-3}))$$ and
$$T(B_n)-T(B_{n-1})=T(G_{n-1})-yT(G_{n-2})+y(T(G_{n-2})-xT(G_{n-3}))$$
$$=T(G_{n-1})-xyT(G_{n-3}).$$
The recurrence formula (\ref{recGn}) provides that
$$T(B_n)=(2+x+y)T(B_{n-1})-(1+x+y+xy)T(B_{n-2})+xyT(B_{n-3})$$
thus
$$\forall n\geq 3,\,\left(\begin{array}{c}
T(B_{n})\\T(B_{n-1})\\T(B_{n-2})\end{array}\right)=\left(\begin{array}{ccc}
2+x+y&-(1+x+y+xy)&xy\\1&0&0\\0&1&0\end{array}\right)\left(\begin{array}{c}T(B_{n-1})\\T(B_{n-2})\\T(B_{n-3})\end{array}\right).$$
The eigenvalues are\\
$\mu_1=\frac{x+y+1+\sqrt{(1+x+y)^2-4xy}}2,\quad
\mu_2=\frac{x+y+1-\sqrt{(1+x+y)^2-4xy}}2,\quad\mu_3=1$.\\ The matrix
defining eigenvectors $P$ and its inverse $P^{-1}$ are given by:
$$P=\left(\begin{array}{ccc}\displaystyle\frac{1}{(\mu_1-1)(\mu_2-1)}&\displaystyle\frac{\mu_1^2}{(\mu_1-1)(\mu_1-\mu_2)}&
\displaystyle\frac{\mu_2^2}{(\mu_2-1)(\mu_2-\mu_1)}\\\displaystyle\frac{1}{(\mu_1-1)(\mu_2-1)}&\displaystyle\frac{\mu_1}{(\mu_1-1)(\mu_1-\mu_2)}&
\displaystyle\frac{\mu_2}{(\mu_2-1)(\mu_2-\mu_1)}\\\displaystyle\frac{1}{(\mu_1-1)(\mu_2-1)}&\displaystyle\frac{1}{(\mu_1-1)(\mu_1-\mu_2)}&
\displaystyle\frac{1}{(\mu_2-1)(\mu_2-\mu_1)}\end{array}\right)$$
and
$$P^{-1}=\left(\begin{array}{ccc}1&-\mu_1-\mu_2&\mu_1\mu_2\\1&-\mu_2-1&\mu_1\mu_2\\1&-\mu_1-1&\mu_1\mu_2
\end{array}\right).$$
Then, it comes
$$\forall n\geq 4,\,\left(\begin{array}{c}T(B_{n})\\T(B_{n-1})\\T(B_{n-2})\end{array}\right)=P\left(\begin{array}{ccc}1&0&0\\0&
\mu_1^{n-3}&0\\0&0&\mu_2^{n-3}\end{array}\right)P^{-1}\left(\begin{array}{c}T(B_{3})\\T(B_{2})\\T(B_{1})\end{array}\right)$$
where $$\left\{\begin{array}{l}
T(B_3)=x^3+y^3+3x^2+3y^2+4xy+2x+2y\\
T(B_2)=x^2+y^2+xy+x+y\\
T(B_1)=xy\end{array}\right.$$
We can conclude that
$$T(B_n)=(xy-x-y-1)+\mu_1^n+\mu_2^n.$$
\noindent On the hyperbola $(x-1)(y-1)=q$, these eigenvalues are of the form introduced before with $a=1$ for $\{\mu_1,\mu_2\}$ and $a=q$ for the eigenvalue 1. The Tutte polynomial can be written as :
 $$T(B_n)=(q-2)[\lambda_q^-(z)]^n+[\lambda_1^+(z)]^n+[\lambda_1^-(z)]^n$$
where $\lambda^{-}_{q}(z)=1$.
\begin{prop}
For the family of graphs $(B_{n})_{n\geq 0}$, the location of the degeneration of the dominant eigenvalue is described in the complex plane using  the variable $z$  as follows:\\
- $\forall q\in[1,5]$, $q\neq 2$, when $|\lambda_1^+(z)|=|\lambda_1^-(z)|>1$ or $|\lambda_1(z)|=1$
$$\Rightarrow \left\{\begin{array}{lll}
c\in[-q,-1+2\sqrt{q-1}]&{\rm and}&d=0\\
{\rm or}&&\\
c\in[-(q+5)/2,-q]&\rm{and}&d^2=-(c+q)^2\frac{2c+q+5}{2c+q+1}\end{array}\right.$$
- $\forall q>5$, when $|\lambda_1^+(z)|=|\lambda_1^-(z)|>1$ $$\Rightarrow c\in[-1-2\sqrt{q-1},-1+2\sqrt{q-1}]\quad{\rm and}\,d=0.$$
- For $q=2$, when $|\lambda_1^+(z)|=|\lambda_1^-(z)|$ $$\Rightarrow c\in[-3,1],\,d=0.$$
\end{prop}
\noindent{\bf Proof:} Let us notice that if we denote
$$\lambda_{1}(z)=1+\sqrt{q-1}\rho e^{i\alpha}$$ with $\rho\geq 1$ then $$|\lambda_1(z)|= 1+2\sqrt{q-1}\rho \cos\alpha + (q-1){\rho}^2.$$ If  $|\lambda_{1}(z)|=1$, it means 
$$cos(\alpha)=-\sqrt{q-1}\rho/2.$$
As $\rho\geq 1$, the last equation has no solution for $q>5$ which explains the three cases.\\
\noindent For the second and third sentences, we apply directly the Lemma 6.2. For the first one, we have two different cases of degeneration. By using the Lemma 6.2,  $|\lambda_1^+(z)|=|\lambda_1^-(z)|$ gives the set $c=[-1-2\sqrt{q-1},-1+2\sqrt{q-1}]$ and $d=0$. But, on the interval $[-1-2\sqrt{q-1},-q]$, the eigenvalue $1$ is dominant and does not match with a degeneration case.\\
In the same way, the Lemma 6.3 leads us to consider the case $|\lambda_1(z)|=|\lambda_q(z)|$. On the set $c\in[-(q+5)/2,-q]$ and $d^2=-(c+q)^2\frac{2c+q+5}{2c+q+1}$, $\lambda_q(z)=1$ and can be dominant with the eigenvalue $\lambda_1(z)$. However, on the set $c\in[-q,-(q+1)/2]$ and $d^2=-(c+q)^2\frac{2c+q+5}{2c+q+1}$, $\lambda_q(z)=z+q+1$ and can be dominant with the eigenvalue $\lambda_1(z)$ but not appears in the Tutte polynomial. In this case, $\lambda_1(z)$ is the one and only dominant eigenvalue. This is not a case of degeneration.

\begin{figure}[ht]
\begin{center}
\begin{tabular}{cc}
\includegraphics[width=4cm,angle=-90]{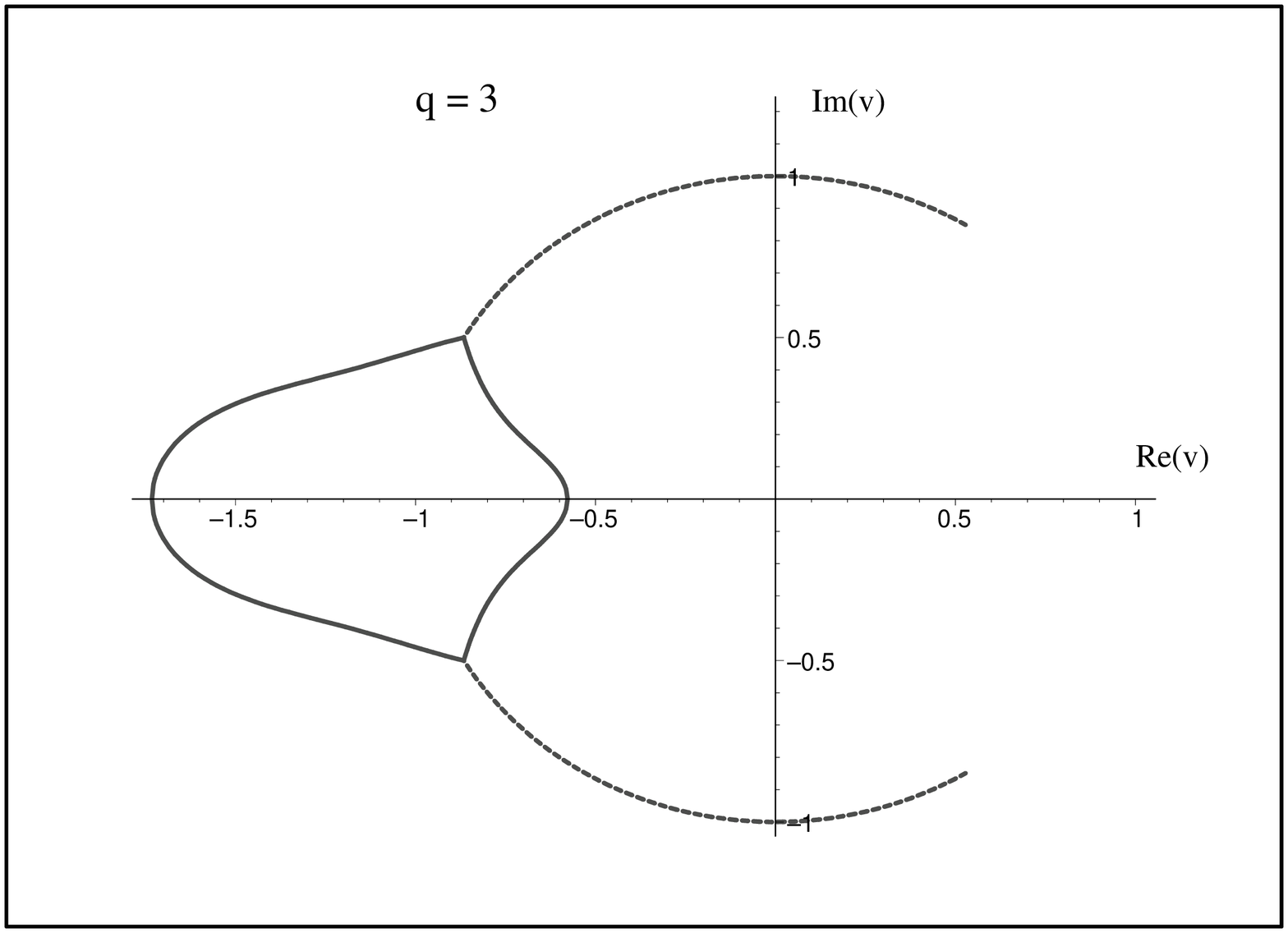}&\includegraphics[width=4cm,angle=-90]{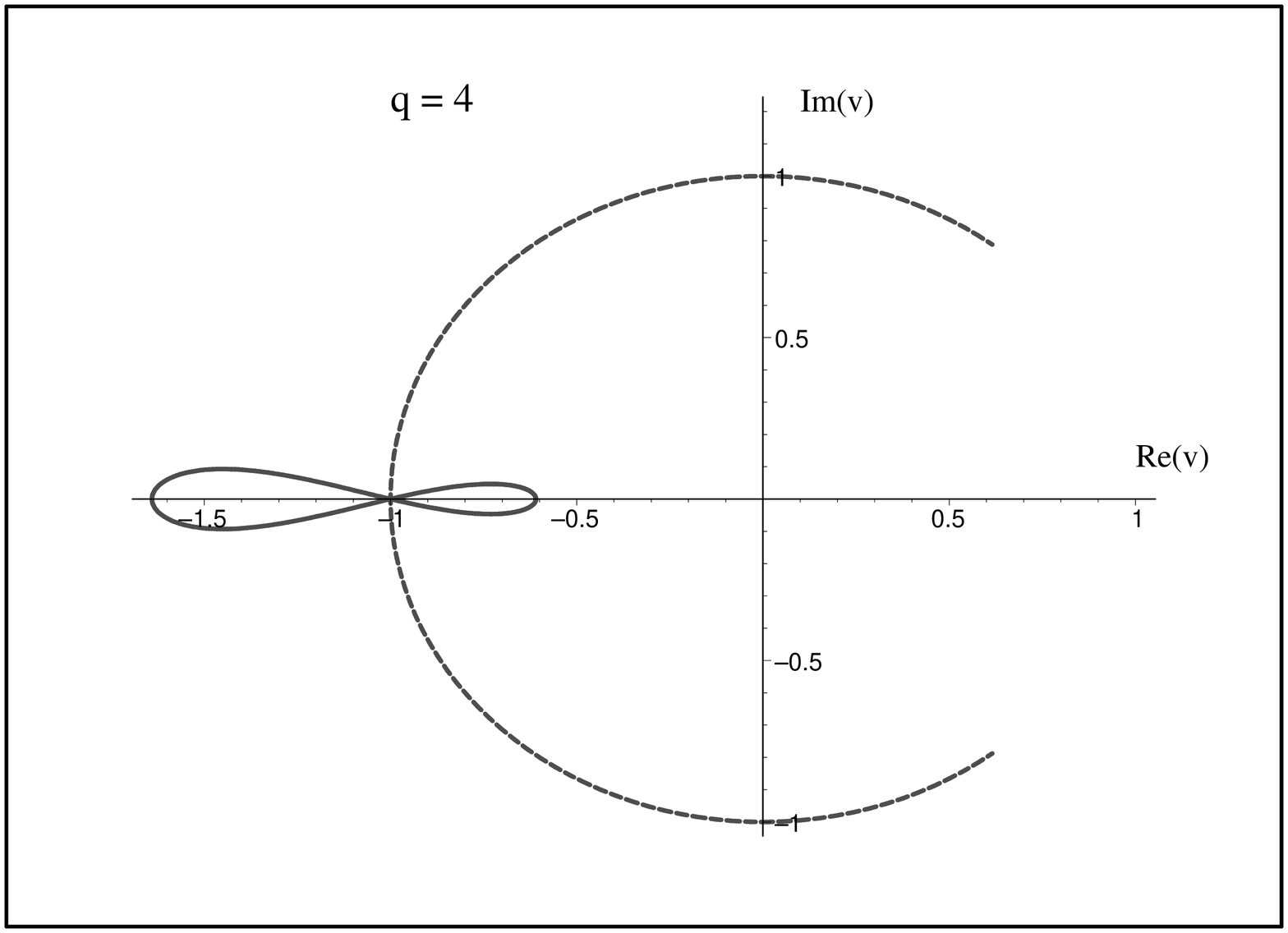}\\
\includegraphics[width=4cm,angle=-90]{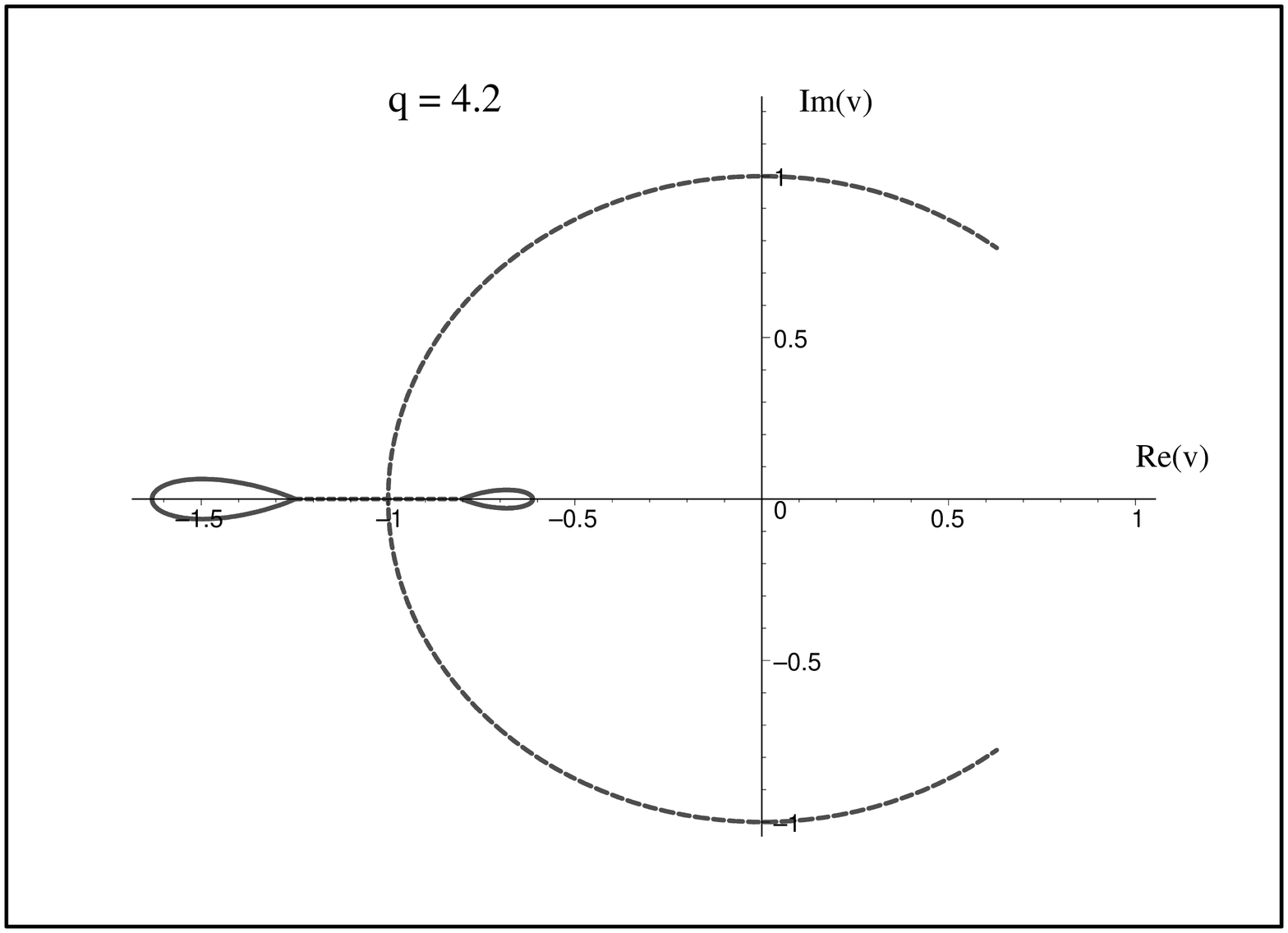}&\includegraphics[width=4cm,angle=-90]{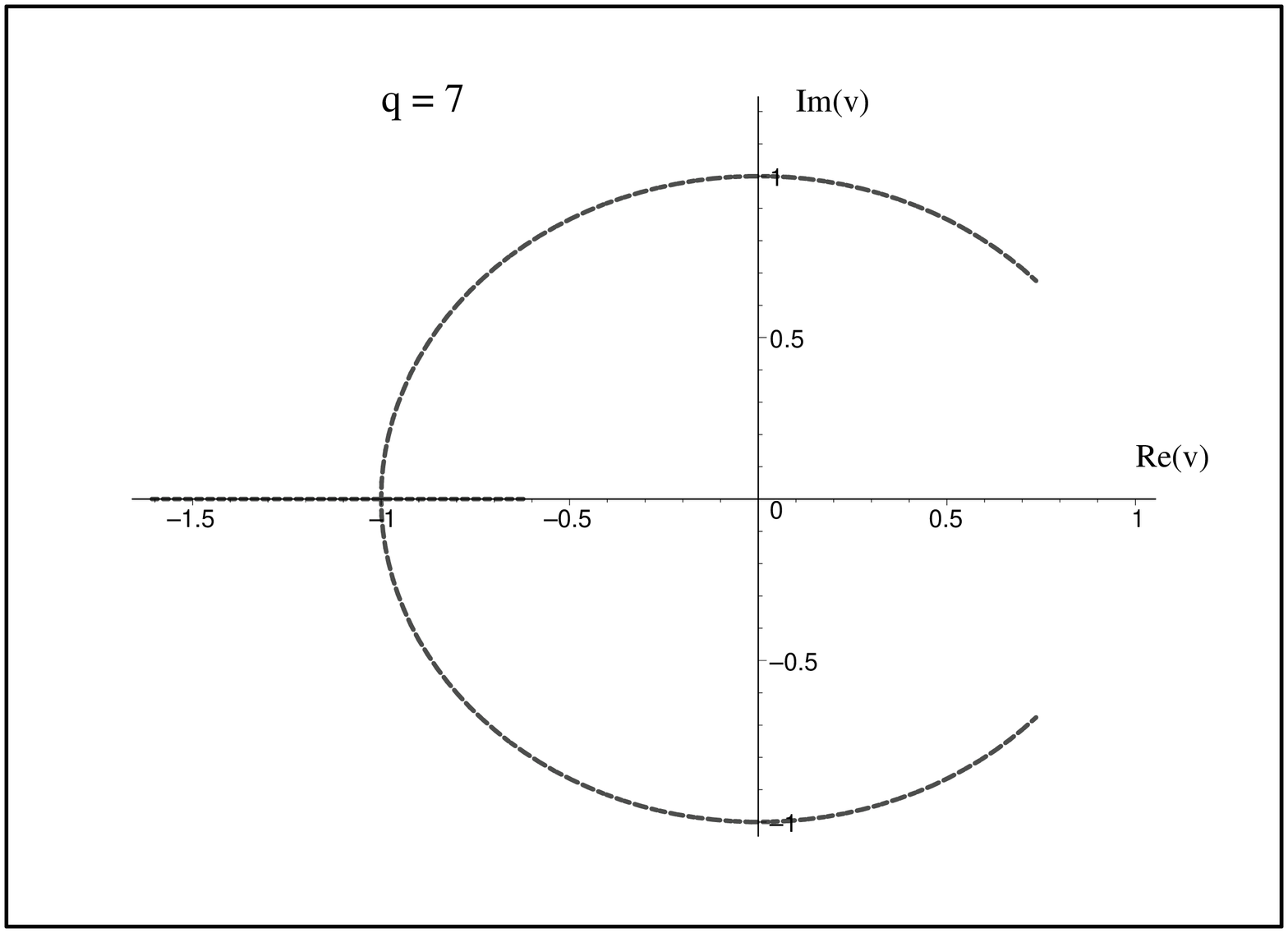}\\
\end{tabular}
\end{center}
\caption{The limiting zero sets for different values of $q$ for the wheel in the $v$ complex plane.} 
\label{wheelq}
\end{figure}

\begin{prop}
For the family of graphs $(B_{n})_{n\geq 0}$, the location of the degeneration of the dominant eigenvalue is described in the complex plane using the variable $v=re^{i\theta}$ as follows:\\
\noindent- if $q\in[1,4]$, $q\neq 2$, $\left\{\begin{array}{l}\theta\in[ \arccos(\frac{-1+2\sqrt{q-1}}{2\sqrt{q}}),\arccos(-\frac{\sqrt{q}}{2})],\,r=1\\
\mbox{or}\,\,
v\in F({\cal P}).\end{array}\right.$\\
\noindent- if $q\in[4,5]$ $\left\{\begin{array}{l}\theta\in[ \arccos(\frac{-1+2\sqrt{q-1}}{2\sqrt{q}}), \pi],\, r=1\\
\mbox{or}\,\,
\theta=\pi,\,r\in[\frac{\sqrt{q}-\sqrt{q-4}}2,\frac{\sqrt{q}+\sqrt{q-4}}2]\\
\mbox{or}\,\,
 v \in F({\cal P}).\end{array}\right.$\\
\noindent- if $q>5$ $\left\{\begin{array}{l}\theta\in[\arccos(\frac{-1+2\sqrt{q-1}}{2\sqrt{q}}), \pi],\,r=1\\
\mbox{or}\,\,
\theta=\pi,\, r\in[1/r_1;r_1]\end{array}\right.$\\
- For $q=2$ $\left\{\begin{array}{l}\theta\in[\arccos(\frac{1}{2\sqrt{2}}), \pi],\,r=1\\
\mbox{or}\,\,
\theta=\pi,\, r\in[\frac{1}{\sqrt{2}};\sqrt{2}]\end{array}\right.$\\
where $r_{1}$, $1/r_1$ are the roots of the polynomial $\sqrt{q}r^{2}-(1+2\sqrt{q-1})r+\sqrt{q}=0$ and ${\cal P}=\{(c,d)\in\R^2,\,c\in[-(q+5)/2,-q],\, d^2=-(c+q)^2(2c+q+5)/(2c+q+1)\}$.
\end{prop}
\noindent{\bf Proof:} We use the proof of the lemma \ref{lem4} to determinate the transformation of the segment $[-q,-1+2\sqrt{q-1}]$ and this one of 
$[-1-2\sqrt{q-1},-1+2\sqrt{q-1}]$. Next, for the set ${\cal P}$, we have only modified the value of $d$ to apply the lemma \ref{lem5}.\\

In the Figure \ref{wheelq}, we show the limiting zero sets for different values of the parameter $q$. The dotted line represents the degeneration case $|\lambda_1^+(z)|=|\lambda_1^-(z)|>1$. The other curve concerns the case where the eigenvalue 1 is dominant. When $q$ is greater than 5, the eigenvalue 1 is never dominant.

\subsection{Cycle with one multiple edge}

\begin{figure}[h!]
\begin{center}
\begin{tabular}{ll}
\includegraphics[width=4cm]{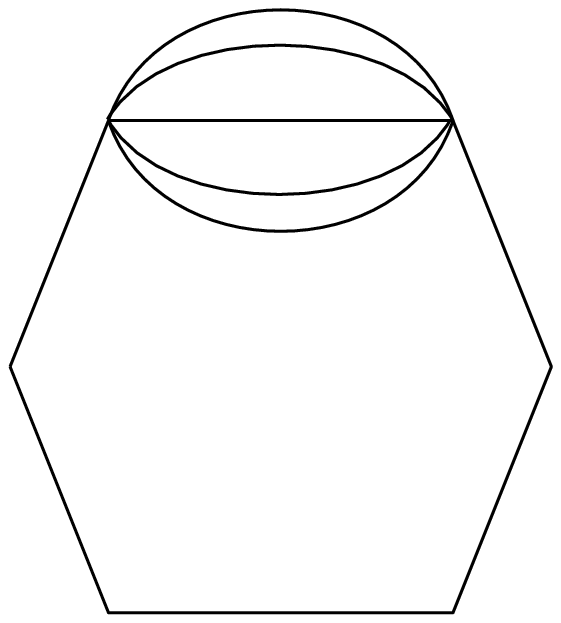}& \includegraphics[width=4cm]{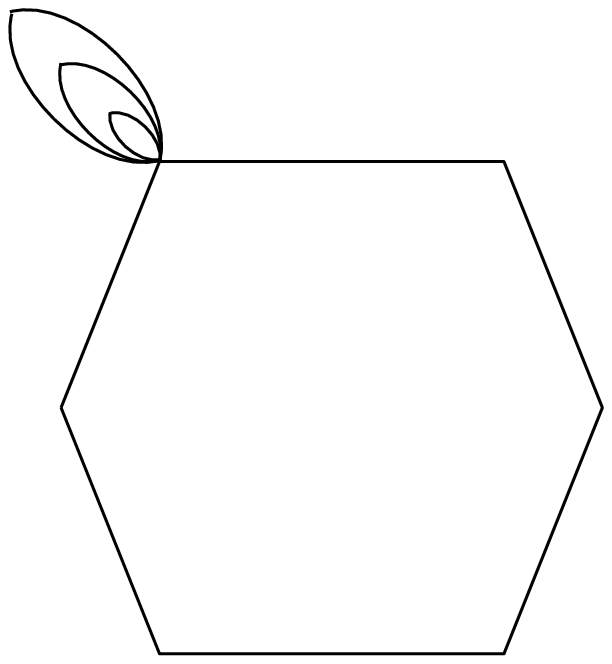} \\
\includegraphics[width=4cm]{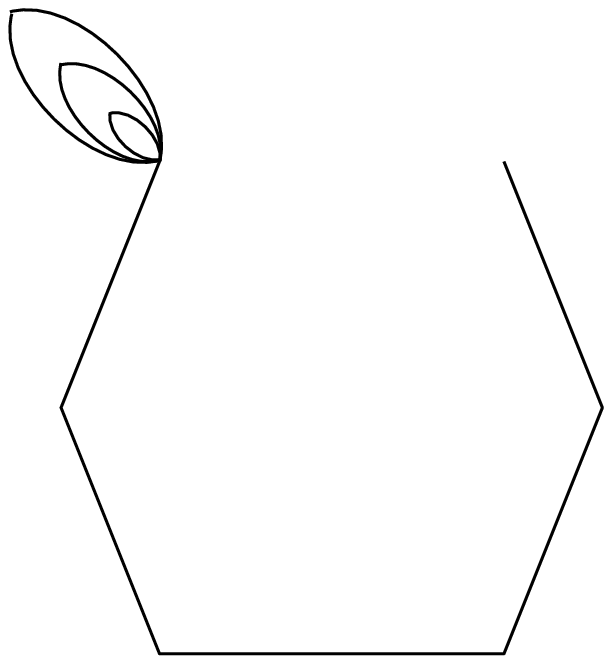} &\includegraphics[width=4cm]{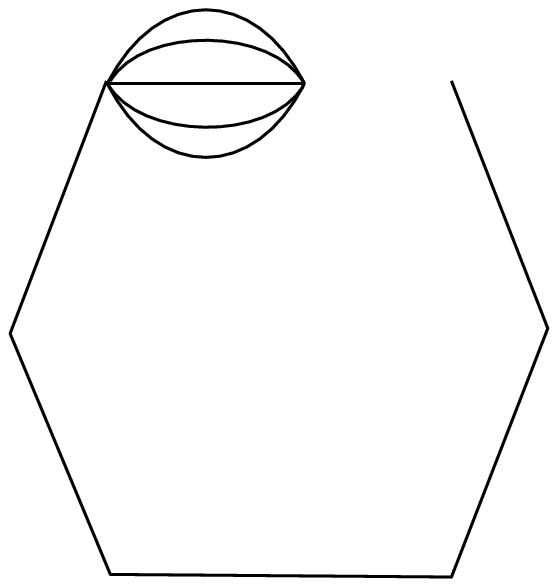}
\end{tabular}
\end{center}\caption{\footnotesize From left to right and from top to bottom, the graphs $G_{n,k}$, $A_{n,k}$, $B_{n,k}$ and $C_{n,k}$ .}
\label{Cycle}
\end{figure}

We consider graphs having $2n$ edges. We assume that, for such a graph, we have one cycle $C_{n+1}$ of length $n+1$ and one edge with order of multiplicity $k=n$. We denote these graphs by $G_{n,k}$.\\

Let us introduce the following graphs $A_{n,k}$, $B_{n,k}$ and $C_{n,k}$: $A_{n,k}$ is the graph having one cycle of length $n$ and one vertex having $k$ loops; $B_{n,k}$ is the graph having one isthmus on length $n$ and ended by $k$ loops; and $C_{n,k}$ is the graph having one isthmus of length $n+1$ and the last edge with order of multiplicity $k$. These graphs are presented in Figure \ref{Cycle}.\\
We have
$$\begin{array}{lll}
T(G_{n,k})&= &T(C_{n-1,k})+T(G_{n-1,k})\\
&=&T(C_{n-1,k-1})+T(B_{n-1,k-1})+T(G_{n-1,k})\\
&=&xT(C_{n-2,k-1})+T(B_{n-1,k-1})+T(G_{n-1,k})\\
&=&x[T(G_{n-1,k-1})-T(G_{n-2,k-1})]+T(B_{n-1,k-1})+T(G_{n-1,k})\\
&=&x[T(G_{n-1,k-1})-T(G_{n-2,k-2}-T(A_{n-2,k-2})]+\\
&&\hspace*{4cm}+\,T(B_{n-1,k-1})\,+\,T(G_{n-1,k})\end{array}$$
then
$$\begin{array}{lll}
T(G_{n,k})&= &x[T(G_{n-1,k-1})-T(G_{n-2,k-2}]-xT(A_{n-2,k-2})\\
&&\hspace*{2cm}+xyT(B_{n-2,k-2})+T(G_{n-1,k-1})+T(A_{n-1,k-1})\end{array}$$
We deduce that
$$\begin{array}{lll}\label{recGnk}
T(G_{n,k})-T(G_{n-1,k-1})&=&x[T(G_{n-1,k-1})-T(G_{n-2,k-2}]\\
&&+(y-x)T(A_{n-2,k-2})+y(x+1)T(B_{n-2,k-2}).
\end{array}$$
Moreover, we have both following relations
\begin{equation}\label{eq1}T(A_{n-1,k-1})=yT(A_{n-2,k-2})+yT(B_{n-2,k-2})\end{equation}
\begin{equation}\label{eq2}T(B_{n-1,k-1})=xyT(B_{n-2,k-2}).\end{equation}
Taking $U_{n,k}=T(G_{n,k})-T(G_{n-1,k-1})$, it comes that \begin{eqnarray}
U_{n,k}=xU_{n-1,k-1}+(y-x)T(A_{n-2,k-2})+y(x+1)T(B_{n-2,k-2})\label{eq3}\\
U_{n-1,k-1}=xU_{n-2,k-2}+(y-x)T(A_{n-3,k-3})+y(x+1)T(B_{n-3,k-3})\label{eq4}\end{eqnarray}
Using the relation (\ref{eq1}) and computing (\ref{eq3})-$y$(\ref{eq4}), we obtain :
$$U_{n,k}-(x+y)U_{n-1,k-1}+xyU_{n-2,k-2}=(xy-1)T(B_{n-2,k-2}).$$
Considering this last relation with order $n-1,k-1$, the equation (\ref{eq2}) and after subtracting and multiplying by $xy$, it leads to
$$U_{n,k}-(x+y+xy)U_{n-1,k-1}+xy(1+x+y)U_{n-2,k-2}-(xy)^2U_{n-3,k-3}=0.$$
Now,  writing with the help of previous Tutte polynomials, we find the following relation:
$$\begin{array}{lll}
T(G_{n,n})&=&\gamma T(G_{n-1,n-1})-(\gamma-1+xy(1+x+y))T(G_{n-2,n-2})\\
&&\hspace*{3cm}+\gamma xyT(G_{n-3,n-3})-(xy)^2T(G_{n-4,n-4})\end{array}$$
where $\gamma=xy+x+y+1$. Hence,
$$\left(\begin{array}{l}
T(G_{n,n})\\
T(G_{n-1,n-1})\\
T(G_{n-2,n-2})\\
T(G_{n-3,n-3})\end{array}\right)=\left(\begin{array}{cccc}\gamma&-(\gamma-1+xy(1+x+y))&\gamma xy&-(xy)^2\\1&0&0&0\\0&1&0&0\\0&0&1&0\end{array}\right)\left(\begin{array}{l}
T(G_{n-1,n-1})\\
T(G_{n-2,n-2})\\
T(G_{n-3,n-3})\\
T(G_{n-4,n-4})\end{array}\right).$$ 

$$\left(\begin{array}{l}
T(G_{n,n})\\
T(G_{n-1,n-1})\\
T(G_{n-2,n-2})\\
T(G_{n-3,n-3})\end{array}\right)=P\left(\begin{array}{cccc}1&0&0&0\\0&x&0&0\\0&0&y&0\\0&0&0&xy\end{array}\right)^{n-3}P^{-1}\left(\begin{array}{l}
T(G_{3,3})\\
T(G_{2,2})\\
T(G_{1,1})\\
T(G_{0,0})\end{array}\right)$$ where
$$T(G_{3,3})=x^3+y^3+xy^{2}+x^2y+x^2y^2+x^2+y^2+xy+x+y,$$
$$T(G_{2,2})=x^2+y^2+xy+x+y, \quad T(G_{1,1})=x+y, \quad T(G_{0,0})=1.$$ 
Then, the Tutte polynomial is given by
$$T(G_{n,n})=\frac{xy-x-y}{(x-1)(y-1)}[x^n+y^n-1]+\frac{(xy)^n}{(x-1)(y-1)}.$$ 

On the hyperbola $(x-1)(y-1)=q$, these eigenvalues are of the form introduced before with $a=0$ for $\{x,y\}$ and with $a=q$ for  $\{1,\,xy\}$. The Tutte polynomial can be written as :
 $$T(G_{n,n})=\frac{1}{q}[\lambda_q^+(z)]^n+\frac{1-q}{q}[\lambda_q^-(z)]^n+\frac{q-1}q[\lambda_0^+(z)]^n+\frac{q-1}q[\lambda_0^-(z)]^n.$$

\begin{prop}
For the family of graphs $(G_{n,n})_{n\geq 0}$, the location of the degeneration of the dominant eigenvalue is described in the complex plane using respectively variable $z$ and $v$ as follows:\\
$\forall q\geq 1$, when $|\lambda_0(z)|=|\lambda_q(z)|$
$$\begin{array}{l}
\Longleftrightarrow\left\{(c,d)\in\R^2,\, c\in[-2-q/2,-q/2]\, , \,d^2= -(2c+q+4)\left[\displaystyle\frac{(c+q)^{2}}{2c+q}\right]\right\}\\
\Longleftrightarrow v\in \Delta_{-\sqrt{q}/2}\cup C((-1/\sqrt{q},0),1/\sqrt{q})\end{array}$$
where $\Delta_{-\sqrt{q}/2}$ is the line $Re(v)=-\sqrt{q}/2$ and $C((-1/\sqrt{q},0),1/\sqrt{q})$ denotes the circle  of center $(-1/\sqrt{q},0)$ and of radius $1/\sqrt{q}$. 
\end{prop}

\begin{pro}
The first part follows again from lemma \ref{lem3} with $a=0$.
Now, using variable $v$, it is easy to see that $|y|=|1+\sqrt{q}v|=1$ leads to $v\in C((-1/\sqrt{q},0),1/\sqrt{q})$ and $|x|=|1+\sqrt{q}/v|=1$ to $v\in \Delta_{-\sqrt{q}/2}$.
\end{pro}

We notice that, for $q>4$ or $-q<-2\sqrt{q}$, then the unit circle $|v|=1$ is not at all in the set of degeneration of the dominant eigenvalue. This is completely the opposite of what we are waiting for in the case of self dual strip of the finite square lattice.

\section{Concluding remarks}
1. Given a Tutte polynomial written with $M$ pairs of eigenvalues like $(\lambda_{a_k}^+(z),\lambda_{a_k}^-(z))$, $k=1\ldots M$. The cases of degeneracy look like to $|\lambda_{a_k}^+(z)|=|\lambda_{a_k}^-(z)|$ or $|\lambda_{a_k}^+(z)|=|\lambda_{a_{k'}}^-(z)|$, for some $k,k'$ in $\{1,\ldots, M\}$, where these eigenvalues are dominant. As the products $|\lambda_{a_k}^+(z)\lambda_{a_k}^-(z)|$ are constant, the first case implies that all eigenvalues have the same magnitude. By using the lemma \ref{lem2}, the degeneracy region in $z$ is given by $\displaystyle\bigcap_{k=1\ldots M}[-a_k-2\sqrt{q-a_k}; -a_k+2\sqrt{q-a_k}]$.\\
For the second case, by using the proof of the lemma \ref{lem1}, $|\lambda_{a}^+(z)|$ is a not increasing function of $a$, $a\in[a_l\,;a_u[$, on the set $D_{]-\infty; -a_u[}$ and  a not decreasing function of $a$, $a\in[a_l\,;a_u[$, on the set $D_{]-a_l; +\infty[}$. Thus, we have to study only one case of degeneracy of these sets : $|\lambda_{a_l}(z)|= |\lambda_{a_u}(z)|$.\\

\noindent 2. To conclude, in the finite case, the conjecture of \cite{Chen96} is true for the strip of triangle with a double edge and the wheel but it is false for the cycle with an edge of high order of multiplicity. But, as a direct consequence of theorem \ref{theo1},
we conclude that, for all graphs studied of this paper, the accumulation set of zeros in positive half plane $Re(v)\geq 0$ is only a sector of the unit circle $|v|=1$. Notice that from corollary \ref{coro1}, we even derive an analyticity result on the positive real axis  as it is foreseeable for some kind of graphs.
But in the general case without knowing other explicit form of eigenvalue, a complete description of the set of degeneration of the dominant eigenvalue in the positive half plane remains open.\\
 
\noindent 3. In the other hand, we found many other self dual graphs belonging to this framework: these graphs can be seen as a transition between the wheel and the cycle with a multiple edge. Moreover, we point out that other choice of parameter $a$ might be interesting. For example consider a family of self dual strip graphs defined as follows. Let the $L_{x}\times L_{y}$ lattice strip have periodic boundary longitudinal or horizontal boundary condition and connect all of the vertices on the upper side of the strip to a single external vertex, while all of the vertices on the lower side of the strip have a free boundary condition. The case $L_{y}=1$ is the wheel graph discussed previously. For this family of graphs, the Tutte polynomial is already known, in particular for the last but one transfer matrices, the eigenvalues obtained for $L_{y}=2$ and $L_{y}=3$ are respectively roots of the following equations
$$\lambda^{4} -(2(x+y)+3)\lambda^3 +(x^{2}+y^{2}+3(x+y)+4xy+1)\lambda^{2}$$
$$-xy(2(x+y)+3)\lambda+(xy)^{2}=0$$
and $$\lambda^{6} -(3(x+y)+5)\lambda^5 +(3(x^{2}+y^{2})+10(x+y)+9xy+6)\lambda^{4}$$
$$-(x^{3}+ y^{3}+5(x^{2}+y^{2})+9xy(x+y)+20xy +6(x+y))\lambda^{3}$$
$$+xy(3(x^{2}+y^{2})+10(x+y)+9xy+6)\lambda^{2}   
 -x^{2}y^{2}(3(x+y)+5)\lambda+(xy)^{3}=0.$$
 For this kind of graphs, the complex-temperature phase diagrams have been already computed for example in \cite{Chang01}.
 For these equations, we find eigenvalues of the form we have studied in this paper.
Denoting by $b_{r}=2+2cos(2\pi/r)$, the Tutte Beraha numbers, we just have to choose values of $a=b_{5}, 3-b_{5}$ for $L_{y}=2$ and $a=b_{7},2-\sqrt{b_{7}},3-b_{7}+\sqrt{b_{7}}$ for $L_{y}=3$. 
Unfortunately, it seems that the eigenvalues of others transfer matrices was not of the desired form. 
It would be interesting to find what choice of parameters $a,\alpha_{a},\beta_{a}$ allows to include other known family of self dual graphs.

\section{Proofs}
In this section, we provide technical lemma on the degeneration of the dominant eigenvalue.
\subsection{Study of the magnitude of the eigenvalues}
For a given $a$, the eigenvalues $\lambda_a^+(z)$ and $\lambda_a^-(z)$ are the solutions of the following equation :
$$X^2-(z+a+2)X+z+q+1=0$$ with $z=c+i\,d$ a complex not real ($d\neq 0$). We will discuss later on the case $d=0$. We find 
$$\lambda_a^{\pm}(z)=\frac12\left((a+c+2+i\,d\pm\displaystyle\sqrt{(z+a+2)^2-4(z+q+1)}\right).$$ By introducing 
$$\begin{array}{l}
A=Re((z+a+2)^2-4(z+q+1))=(a+c)^2-d^2-4(q-a)\\
B=|(z+a+2)^2-4(z+q+1)|=\sqrt{A^2+4d^2(a+c)^2},\end{array}$$
the eigenvalues can be written under this form:
$$\begin{array}{l}
\lambda_a^{\pm}(z)=\frac12\left(a+c+2\pm\sqrt{\frac{B+A}{2}}\right)+\frac{i}2\left(d\pm\sqrt{\frac{2}{B+A}}d(a+c)\right)
\end{array}.$$
Moreover, by using the relation $B^2-A^2=4d^2(a+c)^2$, it comes
$$\left\{\begin{array}{l}
\lambda^{+}_{a}(z)=\frac12\left(a+c+2+(-1)^n\displaystyle\sqrt{\displaystyle\frac{A+B}{2}}\right)+\frac{i}2\left(d+(-1)^{p}\displaystyle\sqrt{\displaystyle\frac{B-A}2}\right)\\
\lambda^{-}_{a}(z)=\frac12\left(a+c+2-(-1)^n\displaystyle\sqrt{\displaystyle\frac{A+B}{2}}\right)+\frac{i}2\left(d-(-1)^{p}\displaystyle\sqrt{\displaystyle\frac{B-A}2}\right)\end{array}\right.$$where
$$n=\left\{\begin{array}{llll}
0&{\rm if}\quad a+c> 0\\
1&{\rm otherwise}\end{array}\right.\quad {\rm and}\quad p=\left\{\begin{array}{llll}
0&{\rm if}\quad d> 0\\
1&{\rm otherwise}\end{array}\right.$$
We consider the following sets 
$$\begin{array}{l}
D^+_{]u,v[}=\{z=c+i\,d,\,(c,d)\in\R^2,\,u<c<v,\,d>0\}\\
D^-_{]u,v[}=\{z=c+i\,d,\,(c,d)\in\R^2,\,u<c<v,\,d<0\}
\end{array}$$
We can state a lemma on the variation of the magnitude of these eigenvalues:

\begin{lem}\label{lem1} 
given a Tutte polynomial written with $M$ pairs of eigenvalues like $(\lambda_{a_k}^+(z),\lambda_{a_k}^-(z))$, $k=1\ldots M$.\\
$\forall a \in [a_l;\,a_u]$, $|\lambda_{\cdot}(z)|=|\lambda_{\cdot}^+(z)|$ is a not decreasing function of the variable $a$ for all $z\in D^{\bullet}_{[-a_{l},+\infty[}$ for $\bullet =+,-$.\newline
$\forall a \in [a_l;\,a_u]$, $|\lambda_{\cdot}(z)|=|\lambda_{\cdot}^+(z)|$ is a not increasing function of $a$ for all $z\in D^{\bullet}_{]-\infty,-a_{u}-2[}$ for $\bullet =+,-$.
\end{lem}
\noindent{\bf Proof:} The magnitude of $\lambda_{a}^{+}(z)$ is given by:
$$4|\lambda_{a}^{+}(z)|^2=\left(a+c+2+(-1)^n\displaystyle\sqrt{\displaystyle\frac{A+B}{2}}\right)^2+\left(d+(-1)^{p}\displaystyle\sqrt{\displaystyle\frac{B-A}2}\right)^2.$$
Taking $U=\sqrt{\frac{B-A}2}$ and $V=\sqrt{\frac{A+B}2}$, it comes that
$$(2|\lambda_{a}^{+}(z)|^2)'=(a+c+2+(-1)^nV)(1+(-1)^nV')+(d+(-1)^{p}U)(-1)^{p}U'$$ then
$$(2|\lambda_{a}^{+}(z)|^2)'=(a+c+2+(-1)^nV)\left(1+(-1)^n\frac{A'+B'}{4V}\right)+(-1)^{p}(d+(-1)^{p}U)\frac{B'-A'}{4U}.$$
Moreover $B'=\frac{1}{B}(AA'+4d^2(a+c))=\frac{W}{B}$. We find that
\begin{eqnarray*}(8|\lambda_{a}^{+}(z)|^2)'=\frac{1}{VB}[a+c+2+(-1)^nV]\left[4VB+(-1)^n(A'B+W)\right]\\+\frac{(-1)^{p}}{UB}[d+(-1)^{p}\,U](W-A'B).
\end{eqnarray*}
Next, we express this differentiation in function of $U$ and $V$ by computing:
$$\left\{\begin{array}{l}
|a+c|=UV/|d|\quad{\rm thus}\quad a+c=(-1)^{n+p}UV/d\\
B=U^2+V^2\\
A=V^2-U^2\\
A'B+W=4V\left(\frac{(-1)^{n+p}UV^2}d+2V+d(-1)^{n+p}U\right)\\
W-A'B=4U\left(-\frac{(-1)^{n+p}U^2V}d-2U+dV(-1)^{n+p}\right).
\end{array}\right.$$
 We conclude that
\begin{eqnarray*}
(2|\lambda_{a}^{+}(z)|^2)'=\frac{(-1)^{n}V}{B}\left\{\left[V\left(1+\frac{U(-1)^{p}}{d}\right)+2\right]^2+\left[U+(-1)^{p}d\right]^2\right\}.\end{eqnarray*}
Hence it is enough to remark that not increasing or not decreasing depends only on the parity of $n$. Now, we identify the eigenvalue between $\lambda_{a}^{+}(z)$ and $\lambda_{a}^{-}(z)$ with greatest modulus which will be denoted by $\lambda_{a}(z)$. We have 
$$|\lambda_{a}^{+}(z)|=|\lambda_{a}^{-}(z)|\Rightarrow \begin{array}{l}
(-1)^n(a+c+2)\sqrt{\displaystyle\frac{A+B}{2}}=-d(-1)^p\sqrt{\displaystyle\frac{B-A}2}\end{array}$$
Thus, by looking at the sign of the members of the previous equalities:
$$|\lambda_{a}^{+}(z)|=|\lambda_{a}^{-}(z)|\Rightarrow \left\{\begin{array}{l}
z\in D^{-}_{]-a-2,-a[}\\ 
z\in D^{+}_{]-a-2,-a[}
\end{array}\right.$$
It comes by the monotonicity of the eigenvalues
$$\lambda_{a}(z)=\left\{\begin{array}{l}\lambda_{a}^{+}(z)\,{\rm on}\,D^{\bullet}_{]-\infty,-a-2[}\\ 
\lambda_{a}^{+}(z)\,{\rm on}\,D^{\bullet}_{]-a,+\infty[}\end{array}\right.$$
Given a Tutte polynomial like $$f_n(z)=\sum_{k=1}^M\Bigl[\alpha_{a_k}(z)[\lambda_{a_k}^{+}(z)]^n+\beta_{a_k}(z)[\lambda_{a_k}^{-}(z)]^n\Bigr ]$$ and $a_l=\min_{k=1..M}{a_k}$, $a_u=\max_{k=1..M}{a_k}$.\newline
 We can assert that, for $a_l\leq a\leq a_u$, $|\lambda_{\cdot}(z)|$ is a not decreasing function of the variable $a$  for all $z\in D^{\bullet}_{]-a_{l},+\infty[}=\displaystyle\bigcap_{k=1\ldots M} D^{\bullet}_{]-a_{k},+\infty[}$. But also that, for $a_l\leq a\leq a_u$, $|\lambda_{\cdot}(z)|$ is a not increasing function of $a$
for all $z\in D^{\bullet}_{]-\infty,-a_{u}-2[}=\displaystyle\bigcap_{k=1\ldots M} D^{\bullet}_{]-\infty,-a_{k}-2[}$ for $\bullet =+,-$. This conclude the proof.\\

\noindent{\bf Remarks:} 1) The case $d=0$ implies that, if $(a+c)^2-4(q-a)> 0$, the eigenvalues are reals and differents. We can not have $|\lambda_{a}^{+}(z)|=|\lambda_{a}^{-}(z)|$. However, if $(a+c)^2-4(q-a)< 0$, the eigenvalues are conjugated and we have $|\lambda_{a}(z)|=|\lambda_{a}^{-}(z)|=|\lambda_{a}^+(z)|$. This gives a case of degeneration and it will be discussed in next lemma.\\

2) We notice that, on the set $D^{\bullet}_{]-a-2,-a[}$, the eigenvalue $\lambda_a(z)$ can be identified to be $\lambda_a^+(z)$ or $\lambda_a^-(z)$. So, we can have a degeneration case like $|\lambda_a^-(z)|=|\lambda_a^+(z)|$ or, when we have $M$ couples of eigenvalues, a degeneration case like $|\lambda_{a_j}^+(z)|=|\lambda_{a_i}^-(z)|$ for $i\neq j$. These cases are studied in the next lemma.

\subsection{Study of some degeneration cases}
At first, we consider a function using one and only one pair of eigenvalues. In this case, it should be interested to study the case of the equality of the magnitudes:
\begin{lem}\label{lem2}
$\forall a\in[0,q-1]$
$$|\lambda^{+}_{a}(z)|=|\lambda^{-}_{a}(z)| \Longleftrightarrow d=0,\, c\in[-a-2\sqrt{q-a},-a+2\sqrt{q-a}]$$
$\forall a\in[q-1,q]$
$$|\lambda^{+}_{a}(z)|=|\lambda^{-}_{a}(z)| \Longleftrightarrow \left\{\begin{array}{l}d=0,\, c\in[-a-2\sqrt{q-a},-a+2\sqrt{q-a}]\\
\mbox{or}\,\,z\in C((-q-1,0),1-q+a)\end{array}\right.$$
where $C((-q-1,0),1-q+a)$ denoted the circle of center $(-q-1,0)$ and of radius $1-q+a$.
\end{lem}
\noindent{\bf Proof:} The equality in magnitude is true only for $z\in D^{\bullet}_{]-a-2,-a[}$ for $\bullet=+,-$.
We have to solve this equation $$(-1)^n(a+c+2)\displaystyle\sqrt{\displaystyle\frac{A+B}{2}}=-d(-1)^p\displaystyle\sqrt{\displaystyle\frac{B-A}2}.$$
It comes
\begin{equation}\label{eqs}
A[(a+c+2)^2+d^2]=B[d^2-(a+c+2)^2].
\end{equation}
As these terms must have the same sign, we find
$$A^2[(a+c+2)^2+d^2]^2=[A^2+4d^2(a+c)^2][d^2-(a+c+2)^2]^2$$
thus
$$d^2[A(a+c+2)-(a+c)(d^2-(a+c+2)^2)][A(a+c+2)+(a+c)(d^2-(a+c+2)^2)]=0$$
It can be written as
$$d^2[A(a+c+2)-(a+c)(d^2-(a+c+2)^2)][(c+q+1)^2+d^2-(q-1-a)^2]=0$$
The case $d=0$ implies that the eigenvalues must be conjugated and gives the segment $c\in[-a-2\sqrt{q-a},-a+2\sqrt{q-a}]$.\\
The third term of the product provides the equation of the circle $C((-q-1,0),|1-q+a|)$. But, as $-a-2<c<-a$, this makes sense only if $1-q+a\geq 0$ thus $q-1\leq a\leq q$. \\
The second term implies that
$$d^2=A(a+c+2)/(a+c)+(a+c+2)^2.$$
If we put it into the equation $(\ref{eqs})$, it follows that
$$A[(a+c+2)^2+d^2]=AB(a+c+2)/(a+c).$$
As $-a-2< c< -a$ in the set $D^{\bullet}_{]-a-2,-a[}$, for this value of $d$, both previous terms have opposite sign, then it is not a solution. That ended the proof of the lemma.\\

\noindent{\bf Remark :} An other way to prove this lemma is to notice that
$$\left\{\begin{array}{l}|\lambda_a^{+}(z)\lambda_a^{-}(z)|=|z+q+1|\\
(\lambda_a^{+}(z)-1)(\lambda_a^{-}(z)-1)=q-a.
\end{array}\right.$$
Thus, it is natural to consider the following form of the eigenvalues:
$$\lambda_a^{+}(z)=1+\sqrt{q-a}\rho e^{i\alpha}\quad \lambda_a^{-}(z)=1+\sqrt{q-a}/\rho e^{-i\alpha}.$$
We obtain
$$\begin{array}{ll}
|\lambda_a^{+}(z)|=|\lambda_a^{-}(z)|&\Leftrightarrow \cos(\alpha)=-\displaystyle\frac{\sqrt{q-a}}2(\rho+\frac{1}{\rho})\,\mbox{or}\, \rho=1.\end{array}$$
The case $\rho=1$ describe the situation where the eigenvalues are conjugated.
This means that $$d=0,\, c\in[-a-2\sqrt{q-a},-a+2\sqrt{q-a}].$$
The other case is possible only when $a\geq q-1$.\\
We deduce that if $\rho\neq 1$, for these value of cosinus, 
$$\left\{\begin{array}{l}
|\lambda_a^{+}(z)|=|\lambda_a^{-}(z)|=\sqrt{1-q+a}\\
|\lambda_a^{+}(z)\lambda_a^{-}(z)|=|z+q+1|\end{array}\right.\Rightarrow|z+q+1|=1-q+a.$$
In this remark, we have just proved the first implication of the previous lemma. It turns out that these expressions of the eigenvalues are useful to study the set of zeros for our family of finite graphs and gives directly an idea of the location of zeros.\\

Now, consider a function using at least one pair of eigenvalues like $(\lambda_{a}^+,\lambda_{a}^-)$ and the pair $(\lambda_{q}^+,\lambda_{q}^-)$. The degeneration can happen when the magnitudes of one eigenvalue of each pair are equal: 
\begin{lem}\label{lem3}
$$|\lambda_a(z)|=|\lambda_q(z)| \Longleftrightarrow  \left\{\begin{array}{l}
c\in[-(q+4+a)/2,-(q+a)/2]\\
d^2=-(c+q)^2\displaystyle\frac{2c+q+4+a}{2c+q+a}\end{array}\right.$$
\end{lem}
\noindent{\bf Proof:} First, we remark that $\lambda_q(z)=1$ or $\lambda_q(z)=z+q+1$. Thus, the relation $|\lambda_a(z)|=|\lambda_q(z)|$ becomes
$$|\lambda_a(z)|=1>|z+q+1|\,{\rm or}\,|\lambda_a(z)|=|z+q+1|>1.$$
In all cases, we have to study $|\lambda_a^+(z)|=1$ or $|\lambda_a^-(z)|=1$. It comes
$$\lambda_a^{+}(z)+\lambda_a^{-}(z)=c+2+a+id\Leftrightarrow \lambda_a^{+}(z)+1+\displaystyle\frac{q-a}{\lambda_a^{+}(z)-1}=c+2+a+id.$$
Taking, for example, $\lambda_a^{+}(z)=e^{i\alpha}$, $$e^{i\alpha}+1 + \displaystyle\frac{q-a}{2(1-\cos\alpha)}(\cos\alpha-1-i\sin\alpha)=c+2+a+id.$$
It gives $$\left\{\begin{array}{l}
c=\cos\alpha-1-q/2-a/2\\
d=\sin\alpha\left[1-\displaystyle\frac{q-a}{2-2\cos\alpha}\right]
\end{array}\right.$$
As $\cos\alpha=c+q/2+a/2+1$, assuming that $-(q+4+a)/2\leq c<-(q+a)/2$, we deduce that 
$$d=\pm\sqrt{1-cos^2\alpha}\left[1+\displaystyle\frac{q-a}{2c+q+a}\right]\Rightarrow d^2 =-(c+q)^2\frac{2c+q+4+a}{2c+q+a}$$
\subsection{Correspondence between $z$ and $v$, case $d=0$}
Many authors work with the variable $v=re^{i\theta}$ such as $z=c+i\,d=\sqrt{q}(v+1/v)$. It leads to $c=\sqrt{q}(r+1/r)\cos(\theta)$ and $d=\sqrt{q}(r-1/r)\sin(\theta)$. We call $F$ the transformation from $z$ to $v$: $$v=F(z)=\frac{z\pm\displaystyle\sqrt{z^2-4q}}{2\sqrt{q}}$$
We have to calculate $F(z)$ when $z\in\{(c,d)\in\R^2,\,d=0\,{\rm and}\,c\in[-a-2\sqrt{q-a},-a+2\sqrt{q-a}]\}$ and when $z$ belongs to the circle $$C((-q-1,0),1-q+a)=\{(c,d)\in\R^2,\,(c+q+1)^2+d^2=(1-q+a)^2\}$$ for $q-1\leq a\leq q$.\\

\begin{lem}\label{lem4}
Let define
$$\left\{\begin{array}{l}
g(a)=-\frac{a}{\sqrt{q}}-2\sqrt{1-\frac{a}{q}}\\
h(a)=-\frac{a}{\sqrt{q}}+2\sqrt{1-\frac{a}{q}}\\
r_1=\frac{-g(a)+\sqrt{g(a)^2-4}}2,\,r_2=\frac{-h(a)+\sqrt{h(a)^2-4}}2\\
\end{array}\right.$$
\begin{itemize}
\item $\forall q\geq 1,\,\forall a\leq 4(\sqrt{q}-1)$,
$$d=0,\, c\in[g(a)\sqrt{q},h(a)\sqrt{q}]\Longleftrightarrow 
\left\{\begin{array}{l}
\theta\in[\arccos(\frac{h(a)}{2}), \pi]\\
r=1\end{array}\right. or 
\left\{\begin{array}{l}
\theta=\pi\\
r\in[1/r_1;r_1]\end{array}\right.$$
\item $\forall 1\leq q\leq 4,\,\forall a> 4(\sqrt{q}-1)$,
$$d=0,\, c\in[g(a)\sqrt{q},h(a)\sqrt{q}]\Longleftrightarrow 
\left\{\begin{array}{l}
\theta\in[\arccos(\frac{h(a)}{2}),\arccos(\frac{g(a)}{2})]\\
r=1\end{array}\right.$$
\item $\forall q\geq 4,\,\forall a> 4(\sqrt{q}-1)$,
$$d=0,\, c\in[g(a)\sqrt{q},h(a)\sqrt{q}]\Longleftrightarrow 
\left\{\begin{array}{l}
\theta=\pi\\
r\in[1/r_2;1/r_1]\cup[r_1;r_2]\end{array}\right.$$
\end{itemize}
\end{lem}

\noindent{\bf Proof:} Taking $d=0=\sqrt{q}(r-1/r)\cos(\theta)$ and $c=\sqrt{q}(r+1/r)\cos(\theta)$, we must have
$$g(a)=-\frac{a}{\sqrt{q}}-2\sqrt{1-\frac{a}{q}}\leq (r+1/r)\cos(\theta)\leq h(a)=-\frac{a}{\sqrt{q}}+2\sqrt{1-\frac{a}{q}}.$$

Remark that $g$ is not increasing on $[0,q-1]$ ($g(0)=-2$, $g(q-1)=-(\sqrt{q}+1/\sqrt{q})$) and not decreasing on $[q-1,q]$ ($g(q)=-\sqrt{q}$).\\
$h$ is not increasing on $[0,q]$ and $h(0)=2$, $h(4[\sqrt{q}-1])=-2$ and $h(q)=-\sqrt{q}$.\\
We first look for the values of $\theta$ when $r=1$. We have to solve
$$g(a)/2\leq \cos(\theta)\leq h(a)/2.$$
\noindent - For $q>4$, we have $g(a)/2\leq -1$. Thus, if $a\leq 4(\sqrt{q}-1)$, $\theta\in[\arccos(\frac{h(a)}{2}), \pi]$ and if $a> 4(\sqrt{q}-1)$, $h(a)/2<-1$ which gives no solution.\\
\noindent - For $q<4$, we have $g(a)/2\geq -1$ and $g(4(\sqrt{q}-1))=-2$. Thus, if $a\geq 4(\sqrt{q}-1)$, $\theta\in[\arccos(\frac{h(a)}{2}), \arccos(\frac{g(a)}{2})]$ and if $a< 4(\sqrt{q}-1)$, $\theta\in[\arccos(\frac{h(a)}{2}), \pi]$.\\

\noindent Now, we look for the values of $r$ when $\theta=\pi$. We have to solve
$$-h(a)\leq (r+1/r)\leq -g(a).$$
\noindent - For $q<4$, $-h(a)<2$. If $a>4(\sqrt{q}-1)$, $-g(a)<2$ which gives no solution. If $a<4(\sqrt{q}-1)$, we obtain $r\in[1/r_1;r_1]$.\\
\noindent - For $q>4$, if $a<4(\sqrt{q}-1)$, $-h(a)<2$ which gives $r\in[1/r_1;r_1]$. And, if $a>4(\sqrt{q}-1)$, $-h(a)>2$ which gives $r\in[1/r_2;1/r_1]\cup[r_1;r_2]$.\\
Remark that the case $\theta=0$ only implies $r=1$.\\

The following lemma gives an analytic expression of $v$ when we apply the transformation $F$ on a circle: 
\begin{lem}\label{lem5}
For all $q-1\leq a\leq q$,\\

$z=c+id\in C((-q-1,0),1-q+a))\Leftrightarrow$\\
$$\left\{\begin{array}{l}d^2=(1-q+a)^2-(c+q+1)^2\\
v=\frac1{2\sqrt{q}}\left[c\pm\sqrt{\frac{A+B}{2}}+i\left(d\pm (-1)^p\sqrt{\frac{B-A}2}\right)\right]

\end{array}\right.$$
with $A=c^2-d^2-4q$ and $B=\sqrt{A^2+4d^2c^2}$.
\end{lem}
\noindent{\bf Proof:} We have just to remark that, for a given $z$, as $$v=F(z)=\frac{z\pm\displaystyle\sqrt{z^2-4q}}{2\sqrt{q}},$$ $v$ can be written under this form :
$$v=\frac{\lambda^{\pm}_{0}(z)-1}{2\sqrt{q}}$$ with
$$\left\{\begin{array}{l}
\lambda^{+}_{0}(z)=\frac12\left(c+2+(-1)^n\displaystyle\sqrt{\displaystyle\frac{A+B}{2}}\right)+\frac{i}2\left(d+(-1)^{p}\displaystyle\sqrt{\displaystyle\frac{B-A}2}\right)\\
\lambda^{-}_{0}(z)=\frac12\left(c+2-(-1)^n\displaystyle\sqrt{\displaystyle\frac{A+B}{2}}\right)+\frac{i}2\left(d-(-1)^{p}\displaystyle\sqrt{\displaystyle\frac{B-A}2}\right)\end{array}\right.$$where
$$n=\left\{\begin{array}{llll}
0&{\rm if}\quad c> 0\\
1&{\rm otherwise}\end{array}\right.\quad {\rm and}\quad p=\left\{\begin{array}{llll}
0&{\rm if}\quad d> 0\\
1&{\rm otherwise}\end{array}\right.$$
 

\bibliographystyle{plain}      
\bibliography{./Biblio/tutte2,./Biblio/eb,./Biblio/horst,./Biblio/biblio,./Biblio/markov,./Biblio/env,./Biblio/bbd,./Biblio/tutte}   

\end{document}